      \documentclass[11pt]{amsart} 
      \usepackage[mathscr]{eucal} 
      \usepackage{amsmath,amsfonts} 
      \parskip\smallskipamount 
          \newtheorem{theorem}{Theorem}[section]

      \newtheorem{corollary}[theorem]{Corollary} 
      \newtheorem{lemma}[theorem]{Lemma}

      \newtheorem{remark}[theorem]{Remark} 
      \makeatletter 
      \@addtoreset{equation}{section} 
      \makeatother 

      \newcommand{\BB}{{\mathbb B}} 
      \newcommand{\CC}{{\mathbb C}}

      \newcommand{\ZZ}{{\mathbb Z}} 
      \newcommand{\DD}{{\mathbb D}} 
       
      \newcommand{\FF}{{\mathbb F}} 
      \newcommand{\TT}{{\mathbb T}}

      \newcommand{\cA}{{\mathcal A}} 
      \newcommand{\cB}{{\mathcal B}} 
       
      \newcommand{\cD}{{\mathcal D}}

      \newcommand{\cG}{{\mathcal G}} 
      \newcommand{\cH}{{\mathcal H}} 
      \newcommand{\cK}{{\mathcal K}} 
       
      \newcommand{\cM}{{\mathcal M}} 
       
      \newcommand{\cP}{{\mathcal P}}

      \newcommand{\cX}{{\mathcal X}}

      \newdimen\expt 
      \def\boxit#1{\setbox0\hbox{$\displaystyle{#1}$} 
            \hbox{\lower.4\expt 
       \hbox{\lower3\expt\hbox{\lower\dp0 
            \hbox{\vbox{\hrule height.4\expt 
       \hbox{\vrule width.4\expt\hskip3\expt 
            \vbox{\vskip3\expt\box0\vskip2\expt}%
       \hskip3\expt\vrule width.4\expt}\hrule height.4\expt}}}}}} 
      
\begin{document} 
       \pagestyle{myheadings} 
      \markboth{ Gelu Popescu}{ Multivariable  Bohr   
      inequalities } 

      \title [  Multivariable Bohr  inequalities ] 
      { Multivariable  Bohr  inequalities} 
        \author{Gelu Popescu} 
      \date{November 8, 2004} 
      \thanks{Research supported in part by an NSF grant} 
      \subjclass[2000]{Primary: 47A20, 47A56;  Secondary:
47A13, 47A63} 
      \keywords{Multivariable operator theory, Bohr's inequality, Holomorphic function, Harmonic function, von Neumann inequality, Poisson transform, Noncommutative disc algebra, Noncommutative  analytic Toeplitz algebra, Fej\' er's inequality} 
        
      \address{Department of Mathematics, The University of Texas 
      at San Antonio \\ San Antonio, TX 78249, USA} 
      \email{\tt gpopescu@math.utsa.edu} 
        
      \begin{abstract} 
        Operator-valued
multivariable Bohr type inequalities are obtained  for: 
\begin{enumerate}
\item[(i)] a class of noncommutative holomorphic functions on the open unit ball of $B(\cH)^n$, generalizing the analytic functions on the open unit disc;
\item[(ii)]
the noncommutative disc algebra $\cA_n$ and the noncommutative analytic Toeplitz algebra $F_n^\infty$;
\item[(iii)] a class of noncommutative selfadjoint harmonic 
functions on the open unit ball of $B(\cH)^n$, generalizing the real-valued harmonic  functions on the  open unit disc;
\item[(iv)] the Cuntz-Toeplitz algebra $C^*(S_1,\ldots, S_n)$, the reduced (resp.~full) group $C^*$-algebra $C_{red}^*(\FF_n)$ (resp.~$C^*(\FF_n)$) of the free group with $n$ generators;
\item[(v)]
a class of analytic functions on the open unit ball of $\CC^n$.
\end{enumerate}

The classical Bohr inequality is shown to be a consequence of Fej\' er's inequality for the coefficients
of positive trigonometric polynomials and Haagerup-de la Harpe inequality
for nilpotent operators. 
Moreover, we provide an inequality which, for analytic polynomials on the open unit disc, is sharper than Bohr's inequality.

 \end{abstract} 
      \maketitle 
      \section*{Introduction}

Let $f(z):=\sum\limits_{k=0}^\infty a_kz^k$ be an  analytic function on the open unit disc 
$\DD:=\{z\in\CC: \ |z|<1\}$ such that  $\|f\|_\infty\leq 1$. Bohr's inequality
\cite{B}
asserts  that 
$$\sum_{k=0}^\infty r^k |a_k|\leq 1\quad \text{  for }\  0\leq r\leq \frac{1} {3}.
$$
   Originally, the inequality was obtained  for
 $0\leq r\leq \frac{1} {6}$. The fact that $\frac{1} {3}$ is the best possible constant was obtained independently by
M. Riesz, Schur, and Weiner. Other proofs were  later  
obtained
 by Sidon \cite{S} and Tomic \cite{T}.
 Dixon \cite{D} used Bohr's inequality in connection with
the long-standing problem of characterizing Banach algebras satisfying  the   von Neumann inequality \cite{vN} (see also \cite{P} and \cite{Pi}). In recent years, multivariable analogues of 
Bohr's inequality were considered by several authors (see \cite{Ai}, \cite{BK},
 \cite{DT}, and \cite{PPoS}).
Paulsen and Singh  \cite{PS} used positivity methods to obtain operator-valued generalizations of Bohr's inequality in the single variable case.

We obtain in this paper 
operator-valued generalizations of Bohr's inequality in
multivariable settings.
      Let $H_n$ be an $n$-dimensional complex  Hilbert space with orthonormal 
      basis 
      $e_1$, $e_2$, $\dots,e_n$, where $n\in\{1,2,\dots\}$.        We consider the full Fock space  of $H_n$ defined by 
      $$F^2(H_n):=\bigoplus_{k\geq 0} H_n^{\otimes k},$$ 
      where $H_n^{\otimes 0}:=\CC 1$ and $H_n^{\otimes k}$ is the (Hilbert) 
      tensor product of $k$ copies of $H_n$. 
      Define the left creation 
      operators $S_i:F^2(H_n)\to F^2(H_n), \  i=1,\dots, n$,  by 
      $$ 
       S_i\varphi:=e_i\otimes\varphi, \quad  \varphi\in F^2(H_n). 
      $$ 
      Let $\FF_n^+$ be the unital free semigroup on $n$ generators 
      $g_1,\dots,g_n$, and the identity $g_0$. 
      The length of $\alpha\in\FF_n^+$ is defined by 
      $|\alpha|:=k$ if $\alpha=g_{i_1}g_{i_2}\cdots g_{i_k}$, and 
      $|\alpha|:=0$ if $\alpha=g_0$. 
       If $T_1,\dots,T_n\in B(\cH)$, the algebra of all bounded operators on a Hilbert space $\cH$,  define 
      $T_\alpha :=  T_{i_1}T_{i_2}\cdots T_{i_k}$ 
      if $\alpha=g_{i_1}g_{i_2}\cdots g_{i_k}$, and 
      $T_{g_0}:=I_\cH$.  

      The    noncommutative analytic Toeplitz algebra   $F_n^\infty$ 
        and  its norm closed version, 
        the noncommutative disc 
       algebra  $\cA_n$,  were introduced by the author   \cite{Po-von}, \cite{Po-funct}, \cite{Po-disc} in 
      connection 
         with a multivariable noncommutative von Neumann inequality. 
      $F_n^\infty$  is the algebra of left multipliers of  
      $F^2(H_n)$  and  can be identified with 
       the 
        weakly closed  (or $w^*$-closed) algebra generated by the left creation 
      operators 
         $S_1,\dots, S_n$  acting on   $F^2(H_n)$, 
          and the identity. 
           The noncommutative disc algebra $\cA_n$ is 
          the  norm closed algebra generated by  
         $S_1,\dots, S_n$,  
          and the identity. 
           When $n=1$, $F_1^\infty$ 
         can be identified 
         with $H^\infty(\DD)$, the algebra of bounded analytic functions 
          on the open unit disc. The  noncommutative analytic Toeplitz algebra $F_n^\infty$ can be viewed as a 
           multivariable noncommutative 
          analogue of $H^\infty(\DD)$. 
       There are many analogies with the invariant 
        subspaces of the unilateral 
       shift on $H^2(\DD)$, inner-outer factorizations, 
        analytic operators, Toeplitz operators, $H^\infty(\DD)$--functional 
         calculus, bounded (resp.~spectral) interpolation, etc. 
      %
%
        
 In Section 1, we show that Bohr's inequality can be improved for analytic polynomials.
More precisely,
we prove that if \ $p(z):=\sum\limits_{k=0}^{m-1} a_k z^k$ is a polynomial with $\|p\|_\infty\leq 1$, then
\begin{equation*}
\sum_{k=1}^{m-1} |a_k| r^k\leq 1\quad 
\text{ for  }\ 0\leq r\leq t_m,
\end{equation*}
 where
 $t_m \in (0,1]$ is the  solution of the equation 
\begin{equation}\label{eq}
\sum_{k=1}^{m-1} t^k\cos\frac{\pi} {\left[\frac{m-1}{k}\right]+2}=\frac{1} {2},
\end{equation}
where $[x]$ is the integer part of $x$.
Moreover, $\{t_m\}_{m=2}^\infty$ is a strictly decreasing sequence which converges to $\frac {1} {3}$.
The above inequality is  a particular case of a more general
multivariable Bohr type inequality, which is  obtained in Section 1,  for noncommutative holomorphic  functions on the open unit unit ball of $[B(\cX)^n]_{<1}$, i.e.,
$$
[B(\cX)^n]_{<1}:=\{(X_1,\ldots, X_n)\in B(\cX)^n:\ \|X_1X_1^*+\cdots +X_nX_n^*\|<1\},
$$
where $\cX$ is an arbitrary Hilbert space. We say that 
$F:[B(\cX)^n]_{<1}\to B(\cX)$ is a universal holomorphic function on $[B(\cX)^n]_{<1}$ 
with scalar coefficients  if there exists a sequence 
$\{a_{\alpha}\}_{\alpha\in \FF_n^+}\subset \CC $ such that
$$
F(X_1,\ldots, X_n)=\sum_{k=0}^\infty \sum_{|\alpha|=k} a_\alpha X_\alpha $$
is convergent in norm for any $(X_1,\ldots, X_n)\in [B(\cX)^n]_{<1}$ and any Hilbert space $\cX$.
We showed in \cite{Po-holomorphic}
that the algebra of all bounded  holomorphic functions on 
$[B(\cX)^n]_{<1}$ can be identified with the noncommutative analytic Toeplitz algebra  $F_n^\infty$, while the subalgebra of all holomorphic functions on 
$[B(\cX)^n]_{<1}$ and continuous on $[B(\cX)^n]_{\leq 1}$
 can be identified with the noncommutative disc algebra $\cA_n$. Many other classical results concerning the analytic functions on the open unit disc  $\DD$ were extended to this noncommutative setting in \cite{Po-holomorphic}.
	
In Section 1, we also consider
multivariable Bohr type inequalities for noncommutative harmonic  functions on the open unit unit ball of $[B(\cX)^n]_{<1}$. We mention that $G$ is a selfadjoint harmonic 
function on $[B(\cX)^n]_{<1}$ if there is a universal holomorphic function $F$  on $[B(\cX)^n]_{<1}$ such that 
$G(X_1,\ldots, X_n)=\text{Re}\, F(X_1,\ldots, X_n)$ for any $(X_1,\ldots, X_n)\in [B(\cX)^n]_{<1}$. 
 As consequences, we obtain 
Bohr type inequalities for the Cuntz-Toeplitz algebra $C^*(S_1,\ldots,S_n)$, 
the reduced (resp.~full) group $C^*$-algebra $C_{red}^*(\FF_n)$ (resp.~$C^*(\FF_n)$) of the free group with $n$ generators. For example, given
$m=2,3,\ldots, \infty$, we show that if $U_1,\ldots, U_n$  are the canonical  unitaries generating 
$C_{red}^*(\FF_n)$ and 
$$
H:=\sum_{1\leq|\alpha|\leq m- 1} \bar a_\alpha U_\alpha^* +a_0 I + 
\sum_{1\leq|\alpha|\leq m-1}  a_\alpha U_\alpha, \quad a_\alpha\in \CC,
$$
 is a selfadjoint element of $C_{red}^*(\FF_n)$ with $\|H\|\leq 1$, then
$$
\left(\sum_{|\alpha|=k}|a_\alpha|^2\right)^{1/2}\leq (1-|a_0|)\cos\frac{\pi} {\left[\frac{m-1}{k}\right]+2}
\quad \text{ for } k=1,\ldots,m-1,
$$
and
$$
\sum_{1\leq |\alpha|\leq m-1} |\bar a_\alpha| r_\alpha +|a_0|  + 
\sum_{1\leq|\alpha|\leq m- 1}  |a_\alpha| r_\alpha
\leq 1
$$
for  any  $r:=(r_1,\ldots, r_n)$ with $r_1\geq 0,\ldots, r_n\geq 0$, and $\|r\|_2\leq t_m$, where   
$t_m \in (0,1]$ is the  solution of the equation \eqref{eq} if $m<\infty$ and $t_\infty=\frac{1} {3}$.

In Section 2, inspired by the work of Paulsen and Singh \cite{PS} in the single variable case, we obtain operator-valued
multivariable 
Bohr type inequalities for   noncommutative holomorphic functions on the open unit ball of $B(\cX)^n$ with coefficients in $B(\cH)$.
As consequences,  we  obtain
 operator-valued
Bohr   inequalities for 
the noncommutative disc algebra $\cA_n$ and the noncommutative analytic Toeplitz algebra $F_n^\infty$.
In particular, we prove that
if $$F(S_1,\ldots, S_n):= \sum_{k=0}^\infty\sum_{|\alpha|=k}S_\alpha\otimes A_{(\alpha)}, \quad A_{(\alpha)}\in B(\cH), $$ is in $F_n^\infty\bar\otimes B(\cH)$, the WOT closed algebra generated by the spatial tensor product,  such that \ $F(0)\geq 0$\  and \ 
Re\,$F(S_1,\ldots, S_n)\leq I$, then
$$
\left\|\sum_{|\alpha|=k}A_{(\alpha)}^*A_{(\alpha)} \right\|^{1/2}\leq 2\|I-A_{(0)}\| \quad \text{ for } k=1,2,\ldots,
$$ and
$$
\sum_{k=0}^\infty\left\|\sum_{|\alpha|=k}T_\alpha\otimes A_{(\alpha)}\right\|\leq
\|A_{(0)}\|+\|I-A_{(0)}\|
$$
for any $n$-tuple of bounded operators $(T_1, \ldots, T_n)\in [B(\cK)^n]_{1/3}$, i.e.,
$$\|T_1T_1^*+\cdots + T_nT_n^*\|^{1/2}\leq \frac {1} {3}.
$$
When $\cH=\CC$ and $A_{(\alpha)}=a_\alpha\in \CC$, we deduce that
$$\sum_{\alpha\in \FF_n^+}  |a_{\alpha}|r_\alpha \leq 1$$
for any $r:=(r_1,\ldots, r_n)$ with $r_1\geq 0,\ldots, r_n\geq 0$, and $\|r\|_2\leq \frac{1} {3}$.
In the single variable case ($n=1$), we find again the classical Bohr inequality \cite{B} and the operator-valued extension of Paulsen and Singh \cite{PS}.
 When $m\geq 2$ and $$F(S_1,\ldots, S_n):= \sum_{k=0}^{m-1}\sum_{|\alpha|=k}S_\alpha\otimes A_{(\alpha)}$$ is a polynomial  such that \ $F(0)\geq 0$\  and \ 
Re\,$F(S_1,\ldots, S_n)\leq I$, then we show that

$$
w\left(A_{(\alpha)}^*:\ |\alpha|=k\right)\leq  2\|I-A_{(0)}\|\cos\frac{\pi} {\left[\frac{m-1}{k}\right]+2}\quad  \text{ for } \ 
1\leq k\leq m-1
$$
and
$$
\sum_{k=0}^{m-1}r^k w\left(A_{(\alpha)}^*:\ |\alpha|=k\right)\leq  \|A_{(0)}\|+\|I-A_{(0)}\|\quad 
\text{ for   }\ 0\leq r\leq t_m,
$$
 where  $t_m \in (0,1]$ is the  solution of the equation 
\eqref{eq} and $w(X_1,\ldots, X_N)$ is the joint numerical radius of $(X_1,\ldots, X_N)$ (see Section 2). We remark that the above  operator-valued Wiener and Bohr type inequalities are new even in the single variable case $n=1$. 

In Section 3, we obtain 
operator-valued Bohr type inequalities for
a class of noncommutative  harmonic 
functions on the unit ball of $B(\cX)^n$ with coefficients in $B(\cH)$.   Wiener and Bohr type inequalities are provided for the coefficients of two harmonic functions on $[B(\cX)^n]_{<1}$ satisfying the inequality 
$H_1(X_1,\ldots, X_n)\leq H_2(X_1,\ldots, X_n)$.
Consequently, we obtain Bohr inequalities  for
  the spatial tensor products 
$C^*(S_1,\ldots, S_n)\otimes B(\cH)$, 
$C_{red}^*(\FF_n)\otimes B(\cH)$, and $C^*(\FF_n)\otimes B(\cH)$.

In Section 4, we provide operator-valued Bohr type inequalities  for
a class of analytic functions on the open unit ball of $\CC^n$. In the scalar case, we obtain Bohr inequalities for the elements of  $H^\infty_{\text{\rm sym}}(\BB_n)$, a Banach  space of analytic functions on $\BB_n$ containing all the polynomials.
More precisely, given $m=2,3,\ldots, \infty$, and
$$f(\lambda_1,\ldots,\lambda_n):=\sum\limits_{{\bf p}\in \ZZ_+^n, |{\bf p}|\leq m-1}
\lambda^{{\bf p}} a_{{\bf p}}, \quad a_{{\bf p}}\in \CC,
$$ 
  an  analytic function in $\BB_n$  such that
 $$\|f\|_{\text{\rm sym}}:=\sup_{0\leq r<1}\left\| f_{\text{\rm sym}}(rS_1,\ldots, rS_n)\right\|<\infty
$$
(see Section 4 for the definition of the symmetrized functional calculus),
 we prove that 
$$\sum\limits_{k=0}^{m-1} \left|  
 \sum\limits_{{\bf p}\in \ZZ_+^n, |{\bf p}|=k} \lambda^{\bf p} |a_{{\bf p}}| \right|\leq  \|f\|_{\text{\rm sym}}
$$
for any $\lambda:=(\lambda_1,\ldots, \lambda_n)\in \BB_n$ with
$\|\lambda\|_2\leq t_m$, where $t_m$ is the solution of the equation \eqref{eq} if $m<\infty$ 
and $t_\infty=\frac{1}{3}$.

\newpage


\section{Bohr inequalities in several variables}

 Haagerup and de la Harpe \cite{HD} proved that any bounded linear
 operator of norm $1$ on a Hilbert space $\cH$ such that $T^m=0$, $m\geq 2$,
 satisfies the inequality $$ \omega(T)\leq \cos \frac {\pi} {m+1}, $$ where
 $\omega(T)$  is the numerical radius of $T$, i.e.,
$$\omega(T):=\{\sup\{|\left< Th,h\right>|:\ h\in \cH, \|h\|=1\}.
$$ 
    They also   showed   that their inequality 
   is equivalent to Fej\' er's inequality 
\cite {Fe} for positive trigonometric polynomials of the form
$$f(e^{i\theta}):= \sum_{k=-m+1}^{m-1} a_k e^{ik\theta}, \quad a_k\in \CC,
$$
which asserts that
$$
|a_1|\leq a_0 \cos \frac {\pi} {m+1}.
$$
\
 In \cite{Po-unitary}, we obtained   multivariable generalizations of the 
 Haagerup--de la Harpe  inequality and   multivariable noncommutative
(resp.~commutative)
 analogues of classical inequalities  (Fej\' er \cite{Fe},
Egerv\' ary-Sz\' azs \cite{ES}) for the coefficients of positive
trigonometric  polynomials. 
  In  particular, we showed that
   if
   $$
   f=\sum_{1\leq|\alpha|\leq m-1} \overline{a}_\alpha S_\alpha^* +a_0 I+
   \sum_{1\leq|\alpha|\leq m-1} a_\alpha S_\alpha
   $$
   is a positive polynomial in $C^*(S_1,\ldots, S_n)$, then
   \begin{equation} \label{fej}
   \left(\sum_{|\alpha|=k} |a_\alpha|^2\right)^{1/2}\leq a_0
   \cos\frac{\pi} {\left[\frac{m-1}{k}\right]+2}
   \end{equation}
    for  $1\leq k\leq m-1$, where $[x]$ is the integer part of $x$.

For $r\geq 0$, we
define 
$$
[B(\cK)^n]_r:=\{(T_1,\ldots, T_n):\ \|T_1T_1^*+\cdots + T_nT_n^*\|^{1/2}\leq r\}.
$$
In what follows  we obtain  a  Bohr type inequality for  analytic polynomials in  the Cuntz-Toeplitz algebra $C^*(S_1,\ldots, S_n)$.

\begin{theorem}\label{sharp} Let $m\geq 2$ and 
let $p(S_1,\ldots, S_n):=   \sum\limits_{|\alpha|\leq m-1} a_\alpha S_\alpha$ be a polynomial such that $p(0)\geq 0$ and 
 $\text{\rm Re}\,p(S_1,\ldots, S_n)\leq I$. Then

\begin{equation*}\sum_{k=0}^{m-1}
\left\| \sum_{|\alpha|=k}|a_\alpha| T_{\alpha}
\right\|\leq 1
\end{equation*}
for any $(T_1,\ldots, T_n)\in [B(\cK)^n]_{t_m}$,
 where
 $t_m \in (0,1]$ is the positive solution of the equation 
\begin{equation}\label{eq1/2}
\sum_{k=1}^{m-1} t^k\cos\frac{\pi} {\left[\frac{m-1}{k}\right]+2}=\frac{1} {2},
\end{equation}
where $[x]$ is the integer part of $x$.
Moreover, $\{t_m\}$ is a strictly decreasing sequence which converges to $\frac {1} {3}$.
In particular,
\begin{equation*}
\sum_{|\alpha|\leq m-1} |a_\alpha| r_\alpha\leq 1
\end{equation*}
for any $r:=(r_1,\ldots, r_n)$ with $r_1\geq 0, \ldots r_n\geq 0$, and $\|r\|_2\leq t_m$.
\end{theorem}

\begin{proof}
The conditions $p(0)\geq 0$ and 
 $\text{\rm Re}\,p(S_1,\ldots, S_n)\leq I$ imply
$$
\sum_{1\leq |\alpha|\leq m-1} -\bar a_\alpha S_\alpha^* +2(1-a_0)I+ \sum_{1\leq |\alpha|\leq m-1} - a_\alpha S_\alpha\geq 0.
$$
According to the inequality \eqref{fej}, we have
\begin{equation}\label{fej2}
 \left(\sum_{|\alpha|=k} |a_\alpha|^2\right)^{1/2}\leq
2(1- a_0)
   \cos\frac{\pi} {\left[\frac{m-1}{k}\right]+2}
   \end{equation}
    for  $1\leq k\leq m-1$.
Let $(T_1,\ldots, T_n)\in [B(\cK)^n]_r$,
 where
 $0\leq r\leq t_m $  and $t_m$ is the positive solution of the equation \eqref{eq1/2}.   Using the noncommutative von Neumann inequality \cite{Po-von} for the row contraction
$[r^{-1}T_1, \ldots, r^{-1} T_n]$  and inequality \eqref{fej2}, we deduce that
\begin{equation*}\begin{split}
\sum_{k=0}^{m-1}
\left\| \sum_{|\alpha|=k}|a_\alpha| T_{\alpha}\right\|&\leq
\sum_{k=0}^{m-1}
\left\| r^k\sum_{|\alpha|=k}|a_\alpha| S_{\alpha}\right\|\\
&= \sum_{k=0}^{m-1} r^k\left( \sum_{|\alpha|=k}|a_\alpha|^2\right)^{1/2}\\
&\leq a_0+2(1-a_0)\sum_{k=1}^{m-1} r^k\cos\frac{\pi} {\left[\frac{m-1}{k}\right]+2}\\
&\leq
a_0+2(1-a_0)\sum_{k=1}^{m-1} t_m^k\cos\frac{\pi} {\left[\frac{m-1}{k}\right]+2}\\
&=a_0+(1-a_0)=1.
 \end{split}
\end{equation*}

For each $m\geq 2$, define the function $f_m:[0,1]\to [0, \infty)$ by setting
$$
f_m(t):= \sum_{k=1}^{m-1} t^k\cos\frac{\pi} {\left[\frac{m-1}{k}\right]+2}.
$$
Notice that $f_m(0)=0$ and $f_m(1)>\cos\frac{\pi}{3}=\frac {1} {2}$.
Since $f_m$ is strictly increasing and continuous, the equation $f_m(t)=\frac {1} {2}$ has a unique solution $t_m\in (0,1]$.
On the other hand, notice that $f_m(t)<f_{m+1} (t)<f(t)$
for any $t\in [0,1)$ and $m\geq 2$, where $f(t)=\sum\limits_{k=1}^\infty t^k$. Since $f(\frac{1}{3})=\frac{1}{2}$, it is clear that $t_m>\frac{1}{3}$ and the sequence $\{t_m\}_{m=2}^\infty$ is strictly decreasing.

A closer look at the  the sequence $\{f_m\}$ reveals that it is uniformly convergent to $f$ on any interval 
$[0,\delta]$ with $0<\delta<1$. Notice also that 
$f_m^\prime (t)\geq \cos\frac{\pi}{m+1}\geq \frac{1}{2}$
for any $m\geq 2$  and $t\in [0,1]$. Applying Lagrange mean value theorem  to the function $f_m$ on the interval $\left[\frac{1}{3}, t_m\right]$, we find $\xi\in\left(\frac {1}{3},t_m\right)$ such that
\begin{equation*}
\begin{split}
\frac{1}{2} -f_m\left(\frac{1}{3}\right)&=
f_m(t_m)-f_m\left(\frac{1}{3}\right)=f^\prime_m(\xi_m)\left(t_m-\frac{1}{3}\right)\\
&\geq \frac{1}{2} \left(t_m-\frac{1}{3}\right)>0
\end{split}
\end{equation*}
Hence, and since $f_m\left(\frac{1}{3}\right)\to f\left(\frac{1}{3}\right)=\frac{1}{2}$, as $m\to\infty$, we deduce that $t_m\to \frac{1}{3}$.
The proof is complete.
 \end{proof}

 Simple computations reveal that        $t_2=1$ and $t_3=
\frac{\sqrt{6}-\sqrt{2}} {2}$.
We  also remark that if $\|p(S_1,\ldots, S_n)\|\leq 1$
then $\text{\rm Re}\,p(S_1,\ldots, S_n)\leq 1$ and,   consequently, Theorem \ref{sharp} holds.

\begin{lemma} \label{poi} Let 
$f(z)=\sum\limits_{k=0}^\infty a_k z^k$ be an analytic function on the open unit disc. Then
$\text{\rm Re}\,f(z)\leq 1$ for any $z\in \DD$ if and only if
\ $\text{\rm Re}\,f(rS)\leq 1$ for $0\leq r<1$, where $S$ is the unilateral shift. Moreover, if $p$ is an analytic  polynomial, then $\text{\rm Re}\,p(z)\leq 1$ if and only if
$\text{ \rm Re}\,p(S)\leq 1$.
\end{lemma}

\begin{proof} Assume that
$\text{ \rm Re}\,f(z)\leq 1$ for any $z\in \DD$.
Notice that, for each $r\in [0,1)$, the function
$f_r(e^{it}):= \sum\limits_{k=0}^\infty r^k a_k e^{ikt}$ 
is in the Hardy space  $H^\infty(\TT)$.
Moreover, $f_r(S):= \sum\limits_{k=0}^\infty r^k a_k S^k$
is convergent in the operator norm.
For every $h(e^{it})\in H^2(\TT)$, we have
\begin{equation*}\begin{split}
\left<[2I-(f_r(S)+f_r(S)^*)]\right.&\left.h(e^{it}),h(e^{it})\right>_{H^2(\TT)}\\
&=
\frac{1} {2\pi} \int_{-\pi}^\pi 
[2-(f_r(e^{it})+\overline{f_r(e^{it})})] |h(e^{it})|^2dt\geq 0.
\end{split}
\end{equation*}
Therefore,  we have $\text{ \rm Re}\,f(rS)\leq 1$ \ for $0\leq r<1$.

Conversely, assume that \ $f_r(S)^*+f_r(S)\leq 2I$  for any
$r\in [0,1)$. Using the Poisson transform of \cite{Po-poisson} associated with the contraction $T:=zI$, where $z\in \DD$, we deduce that
$\overline{f_r(z)}+ f_r(z)\leq 2$ for any  $z\in \DD$ and $r\in [0,1)$. Therefore, $\text{ \rm Re}\,f(z)\leq 1$ for any $z\in \DD$. 
The last part of the lemma is now obvious.
The proof is complete.
\end{proof}

  When $n=1$, Theorem \ref{sharp} and Lemma \ref{poi} imply the following result.

\begin{corollary}  \label{poly}  Let $m\geq 2$ and let   
$p(z)=\sum\limits_{k=0}^{m-1} a_k z^k$ be  an analytic polynomial on the open unit disc such that $p(0)\geq 0$ and 
$\text{ \rm Re}\,p(z)\leq 1$ for $z\in \DD$. Then
$$
|a_k|\leq (1-a_0) \cos\frac{\pi} {\left[\frac{m-1}{k}\right]+2}\quad \text{ for }\ 1\leq k\leq m-1
$$
and 
\begin{equation}\label{b-o}
\sum\limits_{k=0}^\infty |a_k| r^k\leq 1\quad \text{ for } \ 0\leq  r\leq t_m, \end{equation}
 where $t_m$ is the solution of the equation \eqref{eq1/2}.

\end{corollary}

Notice that  if $p $ is a polynomial  with $\|p\|_\infty\leq 1$, then inequality \eqref{b-o} holds and is sharper than Bohr's inequality, when restricted to polynomials.

Now we can prove the following multivariable Bohr type inequality for noncommutative holomorphic functions on the unit open unit ball of $[B(\cX)^n]_{<1}$.

   \begin{theorem}\label{Boh2}
 Let 
       $f(X_1,\ldots, X_n):=\sum\limits_{k=0}^\infty \sum\limits_{|\alpha|=k} 
     a_{\alpha} X_\alpha, 
       \quad a_{\alpha}\in \CC, 
      $ \ 
 be  a  holomorphic function  on  $[B(\cX)^n]_{<1}$    such that  
    $f(0)\geq 0$ and   
$$
\text{ \rm Re}\,f(X_1,\ldots, X_n)\leq I\quad \text{ for any } \ (X_1,\ldots, X_n)\in [B(\cX)^n]_{<1}.
$$        Then 
\begin{equation} \label {win}
\left(\sum_{|\alpha|=k} |a_\alpha|^2\right)^{1/2}\leq
2(1- a_0)\quad \text{ for } \ k=1,2,\ldots,
\end{equation}
 and
      $$\sum_{k=0}^\infty \left\|\sum_{|\alpha|=k} |a_\alpha| 
      T_\alpha\right\|\leq 1 
      $$ 
      for any $(T_1,\ldots, T_n)\in [B(\cK)^n]_{1/3}$. 
In particular, $$\sum\limits_{\alpha\in \FF_n^+} |a_\alpha|r_\alpha\leq  1$$
for any $r:=(r_1,\ldots, r_n)$ with $r_1\geq 0,\ldots, r_n\geq 0$ and $\|r\|_2\leq \frac{1} {3}$.
\end{theorem}
  \begin{proof} 
Since $f$ is holomorphic function on  $[B(\cX)^n]_{<1}$  the series
$\sum\limits_{k=0}^\infty r^k\left(\sum_{|\alpha|=k} |a_\alpha|^2\right)^{1/2}$
is convergent for any $r\in [0,1)$.
Using the noncommutative Poisson transforms of \cite{Po-poisson}, one can easily prove that 
$\text{ \rm Re}\,f(X_1,\ldots, X_n)\leq I$ for any $(X_1,\ldots, X_n)\in [B(\cX)^n]_{<1}$ if and only if  $\text{\rm Re}\,f(rS_1,\ldots, rS_n)\leq I$ for any
 $r\in [0,1)$. Therefore, we have
\begin{equation}\label{susu}
\sum_{k=1}^\infty \sum_{|\alpha|=k} -r^{|\alpha|}\bar a_\alpha S_\alpha^*+2(1-a_0)I + \sum_{k=1}^\infty \sum_{|\alpha|=k} -r^{|\alpha|} a_\alpha S_\alpha
\geq 0.
\end{equation}
For each $m=2,3,\ldots,$ and $r\in [0,1)$, denote 
 $$M_m(r):=\sum\limits_{k\geq m} r^k\left(\sum_{|\alpha|=k} |a_\alpha|^2\right)^{1/2}.$$
 Notice that, for each $r\in [0,1)$, $M_m(r)\to 0$, as $m\to \infty$.
On the other hand, the inequality \eqref{susu} implies
\begin{equation*}
 \sum_{1\leq|\alpha|\leq m-1} -r^{|\alpha|}\bar a_\alpha S_\alpha^*+2(1-a_0+M_m(r))I +  \sum_{1\leq|\alpha|\leq m-1} -r^{|\alpha|} a_\alpha S_\alpha
\geq 0
\end{equation*}
for any $r\in [0,1)$.
According to the inequality \eqref{fej}, we obtain
\begin{equation*}
 r^k\left(\sum_{|\alpha|=k} |a_\alpha|^2\right)^{1/2}\leq
2(1- a_0+ M_m(r))
   \cos\frac{\pi} {\left[\frac{m-1}{k}\right]+2}
   \end{equation*}
    for  $1\leq k\leq m-1$ and any $r\in [0,1)$.
Taking $m\to \infty$ and  then $r\to 1$, we get
the inequality  \eqref{win}.
Hence, and using the noncommutative von Neumann inequality
if 
  $(T_1,\ldots, T_n)\in [B(\cK)^n]_t$ and  $0\leq t\leq \frac{1}{3}$, we obtain 
      \begin{equation*} 
      \begin{split} 
      \sum_{k=0}^\infty \left\| \sum_{|\alpha|=k} |a_{\alpha}| T_\alpha 
      \right\| 
       &\leq 
       \sum_{k=0}^\infty t^k\left\| \sum_{|\alpha|=k} |a_{\alpha}|S_\alpha 
        \right\| 
       \\ 
       &= 
       \sum_{k=0}^\infty t^k\left( \sum_{|\alpha|=k} |a_{\alpha}|^2 
        \right)^{1/2}\\ 
&\leq |a_{0}|+2|1-a_{0}|\frac{t}{1-t}\\
        &\leq |a_{0}|+|1-a_{0}|=1 
      \end{split} 
      \end{equation*} 
      for $0\leq t \leq \frac {1} {3}$. 
This completes the proof.
      \end{proof}

\begin{corollary} \cite{PPoS} \label{Pa-Si}Let 
$f(z)=\sum\limits_{k=0}^\infty a_k z^k$ be an analytic function on the open unit disc such that $f(0)\geq 0$ and 
$\text{ \rm Re}\,f(z)\leq 1$ for $z\in \DD$. Then
$$
\sum\limits_{k=0}^\infty |a_k| r^k\leq 1\quad  \text{ for } \ 0\leq  r\leq \frac {1} {3}.
$$
 \end{corollary}

\begin{proof}
Since $\limsup\limits_{k\to \infty} |a_k|^{1/k}\leq 1$, it is clear that 
$\sum\limits_{n=0}^\infty r^n a_nS^n$ is norm convergent for any $r\in [0,1)$. Now the result follows from 
Lemma \ref{poi} and Theorem \ref{Boh2}.
\end{proof}

Another consequence of Theorem \ref{Boh2} is the following Bohr inequality for the noncommutative analytic Toepltz algebra.

\begin{corollary}\label{Boh3}     If $f(S_1,\ldots, S_n):=\sum\limits_{\alpha\in 
      \FF_n^+} 
       a_\alpha S_\alpha$ is in the noncommutative analytic Toeplitz algebra
$ F_n^\infty$, then 
      $$ \sum_{k=0}^\infty \left\|\sum_{|\alpha|=k} |a_\alpha| 
      T_\alpha\right\|\leq 
      \|f(S_1,\ldots, S_n)\| 
      $$ 
      for any $(T_1,\ldots, T_n)\in [B(\cK)^n]_{1/3}$. 
In particular, $$\sum\limits_{\alpha\in \FF_n^+} |a_\alpha|r_\alpha\leq \|f(S_1,\ldots, S_n)\|$$
for any $r:=(r_1,\ldots, r_n)$ with $r_1\geq 0,\ldots, r_n\geq 0$ and $\|r\|_2\leq \frac{1} {3}$.          %
      \end{corollary} 
    
\begin{proof}
 It follows from  Theorem \ref{Boh2} if we 
      assume 
      that        $\|f(S_1,\ldots, S_n)\|=1$ and notice that, since 
$f(S_1,\ldots, S_n)\in F_n^\infty$, the series
$\sum\limits_{k=0}^\infty \left\|\sum\limits_{|\alpha|=k} r^{|\alpha|}a_\alpha 
      S_\alpha\right\| $ is convergent. 
\end{proof}

      We remark that  Corollary \ref{Boh3} was obtained in 
      \cite{PPoS}, using different methods.   When $n=1$, Corollary \ref{Boh3} and Lemma \ref{poi} imply   the classical Bohr inequality.

\begin{corollary}\label{Bohr's}  $($Bohr's inequality$)$ If $f(z)=\sum\limits_{k=0}^\infty a_k z^k$ is a bounded analytic function on the open unit disc, then
$$
\sum\limits_{k=0}^\infty |a_k| r^k\leq\|f\|_\infty\quad  \text{ for } \ 0\leq  r\leq \frac {1} {3}.
$$
 \end{corollary}

 We remark that, in the particular case $n=1$, the proofs of Theorem \ref{sharp}, Theorem \ref{Boh2},   Corollary \ref{Boh3}, and Corollary \ref{Bohr's}  show that Fej\' er's inequality for the coefficients of positive trigonometric polynomials implies Bohr's inequality for bounded analytic functions on the open unit disc. Here, we should add that  Fej\' er's inequality is equivalent to  Haagerup-de la Harpe inequality, which implies the Egerv\" ary-Sz\' asz inequality (i.e.,
inequality \eqref{fej} in the particular case when $n=1$ and $2\leq k\leq m-1$)
(see \cite{Po-unitary}).

In what follows we obtain Bohr type inequalities for a class of selfadjoint polynomials in the Cuntz-Toeplitz $C^*$-algebra $C^*(S_1,\ldots, S_n)$.

         \begin{theorem} \label{SUU1} 
Let $m\geq2$ and let 
         $$ 
         H(S_1,\ldots, S_n):=\sum_{k=1}^{m-1}\sum_{|\alpha|=k} \overline{a}_\alpha 
      S_\alpha^* 
         +a_0 I+ 
          \sum_{k=1}^{m-1}\sum_{|\alpha|=k} a_\alpha S_\alpha, 
      \quad a_\alpha\in \CC, 
         $$ 
         be a selfadjoint 
element of $C^*(S_1,\ldots, S_n)$
such that          $\|H(S_1,\ldots, S_n)\|\leq 1.$

            Then the following statements hold.
           \begin{enumerate} 
           \item[(i)] 
           $\left(\sum\limits_{|\alpha|=k} |a_\alpha|^2\right)^{1/2}\leq (1-|a_0|)\cos\frac{\pi}{\left[\frac{m-1}{k}\right]}$ \ for  \ 
      $1\leq k\leq m-1$, 
where $[x]$ is the integer part of $x$.
           \item[(ii)] 
           $\sum\limits_{k=0}^{m-1} \left\|\sum\limits_{|\alpha|=k} 
           |a_\alpha| T_\alpha \right\|\leq 1$ 
           \ for any \ $[T_1,\ldots, T_n]\in [B(\cH)^n]_{\gamma_m}$,
where
 $\gamma_m \in (0,1]$ is the  solution of the equation 
\begin{equation}\label{eq2}
\sum_{k=1}^{m-1} t^k\cos\frac{\pi} {\left[\frac{m-1}{k}\right]+2}=1.
\end{equation}
Moreover, $\{\gamma_m\}$ is a strictly decreasing sequence which converges to $\frac {1} {2}$.
\item[(iii)] 
$\sum\limits_{k=1}^{m-1} \left\|\sum\limits_{|\alpha|=k} 
           |\overline{a}_\alpha| T^*_\alpha \right\|
+
|a_0|+\sum\limits_{k=1}^{m-1} \left\|\sum\limits_{|\alpha|=k} 
           |a_\alpha| T_\alpha \right\|
\leq 1$
\ for any \ $[T_1,\ldots, T_n]\in [B(\cH)^n]_{t_m}$,
where
 $t_m \in (0,1]$ is the  solution of the equation 
\begin{equation}\label{eq3}
\sum_{k=1}^{m-1} t^k\cos\frac{\pi} {\left[\frac{m-1}{k}\right]+2}=\frac{1} {2}.
\end{equation}
 Moreover, $\{t_m\}$ is a strictly decreasing sequence which converges to $\frac {1} {3}$.
           \end{enumerate} 
           In particular, if $r_1,\ldots, r_n\geq 0$, then 
\begin{enumerate}
\item[(iv)]
$\sum\limits_{k=0}^{m-1} \left(\sum\limits_{|\alpha|=k} |a_\alpha| 
      r_\alpha\right)\leq 1 
            $
\ if \ 
            $\|(r_1,\ldots, r_n)\|_2\leq \gamma_m$;   
           \item[(v)]
 $|a_0|+ 2\sum\limits_{k=1}^{m-1} \left(\sum\limits_{|\alpha|=k} |a_\alpha| 
      r_\alpha\right)\leq 1 
$
\ if \ 
            $\|(r_1,\ldots, r_n)\|_2\leq t_m$.
\end{enumerate}
         \end{theorem} 
\begin{proof}
Notice that we can assume that $a_0>0$ and 
$$\sum_{k=1}^{m-1}
\sum_{|\alpha|=k} -\bar a_\alpha S_\alpha^*+ (1-a_0)I +
\sum_{|\alpha|=k} - a_\alpha S_\alpha\geq 0.
$$
Using inequality \eqref{fej}, we deduce (i).
The rest of the proof is similar to that of Theorem
\ref{sharp}. We shall omit it.
 \end{proof}

      Let  $\FF_n$ be the free group with generators $g_1,\ldots, g_n$, and let 
      $\ell^2(\FF_n)$  be the 
      Hilbert space defined by 
      $$ 
      \ell^2(\FF_n):=\{f:\FF_n\to \CC:\sum\limits_{\sigma\in\FF_n} 
      |f(\sigma)|^2<\infty\}. 
      $$ 
      The canonical basis of 
      $\ell^2(\FF_n)$ is $\{\xi_ \sigma\}_{ \sigma\in\FF_n}$, where $\xi_ 
      \sigma(t)=1$ 
      if $t= \sigma$ and $\xi_ \sigma(t)=0$ 
      otherwise. 
      For each $i=1,\dots,n$, let  $U_i\in B(\ell^2(\FF_n))$ be the 
      unitary operator defined by 
      $$ 
      U_i\left(\sum\limits_{ \sigma\in\FF_n}a_ \sigma \xi_ \sigma\right): 
      =\sum\limits_{\sigma\in\FF_n} a_ \sigma \xi_{g_i \sigma}, \qquad 
      \left(\sum\limits_{\sigma\in\FF_n}| a_ \sigma|^2<\infty\right). 
      $$ 
       The reduced  group $C^*$-algebra  $C_{\text{\rm red}}^*(\FF_n)$ 
      is the $C^*$-algebra generated by $U_1,\dots, U_n$.

\begin{corollary}\label{CUU}
Let 
         $$ 
         H(U_1,\ldots, U_n):=\sum_{k=1}^{m-1}\sum_{|\alpha|=k} \overline{a}_\alpha 
      U_\alpha^* 
         +a_0 I+ 
          \sum_{k=1}^{m-1}\sum_{|\alpha|=k} a_\alpha U_\alpha, 
      \quad a_\alpha\in \CC, 
         $$ 
         be a selfadjoint element of 
           $C_{\text{\rm red}}^*(\FF_n)$ with           $\|H(U_1,\ldots, U_n)\|\leq 1.$ 
Then the conclusion of Theorem $\ref{SUU1}$ holds.
\end{corollary} 
\begin{proof}

 Notice that 
      the Hilbert space $\ell^2(\FF_n^+)$ can be seen as a subspace of 
      $\ell^2(\FF_n)$ and  the full Fock space $F^2(H_n)$ can be naturally 
      identified 
      to $\ell^2(\FF_n^+)$. 
      Under this identification, we have  
      $U_i|_{F^2(H_n)}=S_i$, $i=1,\dots,n$, 
      where $S_1,\dots,S_n$ are the left creation operators.  Consequently, we 
      have 
      $$ 
      H(S_1,\ldots, S_n)= P_{\ell^2(\FF_n^+)}H(U_1,\ldots, U_n)| 
      \ell^2(\FF_n^+). 
      $$ 
      Since $\|H(U_1,\ldots, U_n)\|\leq 1$, we have     $\|H(S_1,\ldots, S_n)\|\leq 1$.
Applying now Theorem \ref{SUU1}, the result follows. 
\end{proof}
 
  We recall that the full group $C^*$-algebra $C^*(\FF_n)$ is generated by an $n$-tuple of 
      universal 
      unitaries ${\bf U}_1,\ldots, {\bf U}_n$. 
We remark that a result similar to Corollary \ref{CUU} holds for $C^*(\FF_n)$. Indeed, assume that 
$\|H({\bf U}_1,\ldots, {\bf U}_n)\|\leq 1$.
        Due to the universal property of  the unitaries       ${\bf U}_1,\ldots, {\bf U}_n$, there is a $*$-representation $\pi$ of 
       $C^*({\bf U}_1,\ldots, {\bf U}_n)$ onto $C^*(U_1,\ldots, U_n)$ such 
       that $\pi({{\bf U}}_i)=U_i$, 
       \ $i=1,\ldots, n$. Therefore, 
       $ 
       H(U_1,\ldots, U_n)=\pi(H({\bf U}_1,\ldots, {\bf U}_n))\leq I. 
       $ 
        Applying now Corollary \ref{CUU}, the result follows.

The following result is a multivariable Bohr inequality for selfadjoint harmonic functions on the open unit ball of $[B(\cX)^n]_{<1}$.

         \begin{theorem} \label{SUU} 
Let 
         $$ 
         H(X_1,\ldots, X_n):=\sum_{k=1}^\infty\sum_{|\alpha|=k} \overline{a}_\alpha 
      X_\alpha^* 
         +a_0 I+ 
          \sum_{k=1}^\infty\sum_{|\alpha|=k} a_\alpha X_\alpha, 
      \quad a_\alpha\in \CC, 
         $$ 
         be a  selfadjoint  harmonic function on $[B(\cX)^n]_{<1}$
such that          $\|H(X_1,\ldots, X_n)\|\leq 1$ for any $(X_1,\ldots, X_n)\in [B(\cX)^n]_{<1}$.

            Then 
           \begin{enumerate} 
           \item[(i)] 
           $\left(\sum\limits_{|\alpha|=k} |a_\alpha|^2\right)^{1/2}\leq 1-|a_0|$ \ for any \ 
      $k=1,2,\ldots;$ 
           \item[(ii)] 
           $\sum\limits_{k=0}^\infty \left\|\sum\limits_{|\alpha|=k} 
           |a_\alpha| T_\alpha \right\|\leq 1$ 
           \ for any \ $[T_1,\ldots, T_n]\in [B(\cH)^n]_{1/2};$ 
\item[(iii)] 
$\sum\limits_{k=1}^\infty \left\|\sum\limits_{|\alpha|=k} 
           |\overline{a}_\alpha| T^*_\alpha \right\|
+
|a_0|+\sum\limits_{k=1}^\infty \left\|\sum\limits_{|\alpha|=k} 
           |a_\alpha| T_\alpha \right\|
\leq 1$
\ for any \ $[T_1,\ldots, T_n]\in [B(\cH)^n]_{1/3}$.
           \end{enumerate} 
           In particular, if $r_1,\ldots, r_n\geq 0$, then 
\begin{enumerate}
\item[(iv)]
$\sum\limits_{k=0}^\infty \left(\sum\limits_{|\alpha|=k} |a_\alpha| 
      r_\alpha\right)\leq 1 
            $
\ if \ 
            $\|(r_1,\ldots, r_n)\|_2\leq \frac {1} {2}$;   
           \item[(v)]
 $|a_0|+ 2\sum\limits_{k=1}^\infty \left(\sum\limits_{|\alpha|=k} |a_\alpha| 
      r_\alpha\right)\leq 1 
$
\ if \ 
            $\|(r_1,\ldots, r_n)\|_2\leq \frac {1} {3}$.
\end{enumerate}
         \end{theorem} 
      \begin{proof} 
Since $H(X_1,\ldots, X_n)$ is selfadjoint and $\|H(X_1,\ldots, X_n)\|\leq 1$, we can assume that $a_0\geq 0$.
Therefore, we have
$$ \sum_{k=1}^\infty\sum_{|\alpha|=k}r^{|\alpha|}\overline{b}_\alpha 
      S_\alpha^* 
         +b_0 I+ 
         \sum_{k=1}^\infty\sum_{|\alpha|=k}r^{|\alpha|} b_\alpha S_\alpha \geq 0, \quad r\in [0,1),
$$
where $b_0:=1-a_0$ and $b_\alpha:=-a_\alpha$ if $\alpha\in \FF_n^+$, $|\alpha|\geq 1$. 
As in the proof of Theorem \ref{Boh2}, we deduce the inequality (i).
Hence, it is easy to see that if $r\leq \frac{1}{2}$, then
\begin{equation}\label{a0}
\sum_{k=0}^\infty r^k\left( \sum_{|\alpha|=k} |a_\alpha|^2\right)^{1/2}\leq a_0+\frac{r} {1-r}(1-a_0)\leq 1.
\end{equation}
Let $(T_1,\ldots, T_n)\in [B(\cH)^n]_r$  and assume that $0<r\leq \frac {1}{2}$.  
 Using the noncommutative von Neumann inequality  \cite{Po-von} for the row contraction  $[\frac{1}{r}T_1,\ldots,\frac{1}{r}T_n]$
   and  inequality \eqref{a0}, we get
\begin{equation*}
\begin{split}
\sum_{k=0}^\infty
\left\|\sum\limits_{|\alpha|=k} 
           |a_\alpha| T_\alpha \right\|&\leq 
\sum_{k=0}^\infty r^k\left\|\sum\limits_{|\alpha|=k} 
           |a_\alpha| S_\alpha \right\| \\
&=\sum_{k=0}^\infty r^k\left(\sum\limits_{|\alpha|=k} 
           |a_\alpha|^2 \right)^{1/2}\leq 1,
\end{split}
\end{equation*}
which proves the inequality (ii). To prove (iii), notice that 
\begin{equation*} 
\sum_{k=1}^\infty r^k\left( \sum_{|\alpha|=k} |\overline{a}_\alpha|^2\right)^{1/2}+ \sum_{k=0}^\infty r^k\left( \sum_{|\alpha|=k} |a_\alpha|^2\right)^{1/2}\leq a_0+\frac{2r} {1-r}(1-a_0)\leq 1
\end{equation*}
for  $0\leq r\leq \frac{1}{3}$.  Now, the proof  follows the same lines as the proof of (ii).
The proof is complete.
   \end{proof}

\begin{corollary}\label{CUU3}
Let 
         $$ 
         H(U_1,\ldots, U_n):=\sum_{k=1}^\infty\sum_{|\alpha|=k} \overline{a}_\alpha 
      U_\alpha^* 
         +a_0 I+ 
          \sum_{k=1}^\infty\sum_{|\alpha|=k} a_\alpha U_\alpha, 
      \quad a_\alpha\in \CC, 
         $$ 
         be a selfadjoint element of 
           $C_{\text{\rm red}}^*(\FF_n)$ with           $\|H(U_1,\ldots, U_n)\|\leq 1.$ 
Then the conclusion of Theorem $\ref{SUU}$ holds.
\end{corollary} 

We remark that a result similar to Corollary \ref{CUU3} holds for $C^*(\FF_n)$.

      \bigskip
 
      \section{Operator-valued Bohr  inequalities in several variables}

In this section, we obtain operator-valued
multivariable 
Bohr type inequalities for   noncommutative holomorphic functions on the open unit ball of $B(\cX)^n$,
the noncommutative disc algebra $\cA_n$, and the noncommutative analytic Toeplitz algebra $F_n^\infty$.

We say that a power series

$$
F(S_1,\ldots, S_n):=\sum_{k=0}^\infty \sum_{|\alpha|=k} S_\alpha\otimes A_{(\alpha)}, \quad 
\{A_{(\alpha)}\}_{\alpha\in \FF_n^+}\subset B(\cH), 
$$
generates  a (universal) holomorphic function on 
$[B(\cX)^n]_{<1}$ with coefficients in $B(\cH)$ if 
$$
F(X_1,\ldots, X_n)=\sum_{k=0}^\infty \sum_{|\alpha|=k} X_\alpha
\otimes A_{(\alpha)} $$
is convergent in the operator norm for any $(X_1,\ldots, X_n)\in [B(\cX)^n]_{<1}$ and any Hilbert space $\cX$.
We proved in \cite{Po-holomorphic} that the following statements are equivalent:
\begin{enumerate}
\item[(i)] $F(S_1,\ldots, S_n)$ generates  a   holomorphic function on 
$[B(\cX)^n]_{<1}$;
\item[(ii)]
the series 
$\sum\limits_{k=0}^\infty\left\| \sum\limits_{|\alpha|=k} r^{|\alpha|} 
S_\alpha\otimes A_{(\alpha)} \right\|$ is convergent for any 
$r\in [0,1)$, where $S_1,\ldots, S_n$ are the left creation operators on the full Fock space $F^2(H_n)$;
\item[(iii)]
$\limsup\limits_{k\to\infty} \left\|  \sum\limits_{|\alpha|=k}A_{(\alpha)}^*A_{(\alpha)}\right\|^{1/2k}\leq 1$.
\end{enumerate}

Throughout this section we consider (universal) holomorphic functions on the open unit ball
$[B(\cX)^n]_{<1}$ with coefficients in $B(\cH)$. We also set $A_{(0)}:=A_{(g_0)}$.

First, we  recall a very well-known characterization  of  contractions on Hilbert spaces.

      \begin{lemma}\label{standard} 
      A bounded operator $A:\cK\to \cH$ is a contraction if and only if 
      $$ 
      \left(\begin{matrix} 
      I_\cH & A\\ 
      A^* & I_\cK 
      \end{matrix}\right) 
      $$ 
      is a positive operator acting on the Hilbert space $\cH\otimes \cK$. 
      \end{lemma}

The next positivity  result will be useful in what follows.
      \begin{lemma} \label{posi} 
      Let $P, X_1,\ldots, X_m\in B(\cH)$ and let $V_1,\ldots, V_m\in B(\cK)$ be 
      any 
      isometries with orthogonal ranges. 
      Then 
      \begin{equation} 
      \label{pos1} M(P,X_i):= 
      \left( \begin{matrix} 
      P & X_1^* & X_2^* &\cdots & X_m^*\\ 
      X_1& P &0 &\cdots & 0\\ 
      X_2 &  0& P& \cdots & 0\\ 
      \vdots & \vdots & \vdots& \ddots&  \vdots\\
     X_m&0&0&\cdots&0
      \end{matrix} \right) \geq 0 
      \end{equation} 
      if and only if 
      \begin{equation} 
      \label{pos2}N(P,X_i):= 
      \left( \begin{matrix} 
      I_{\cK}\otimes P & \sum\limits_{i=1}^m V_i \otimes X_i\\ 
      \sum\limits_{i=1}^m V_i^* \otimes X_i^* & I_{\cK}\otimes P 
      \end{matrix} \right) \geq 0. 
      \end{equation} 
      \end{lemma} 
      \begin{proof} 
      Notice that $M(P,X_i)\geq 0$ if and only if $P\geq 0$ and 
      $$ 
      M(I_\cH, (P+\epsilon I_\cH)^{-1/2} X_i (P+\epsilon I_\cH)^{-1/2})\geq 0 
      $$ 
      for any $\epsilon >0$. 
      Applying Lemma \ref{standard} to the row operator 
      $A_\epsilon:=[A_{\epsilon, 1}, \ldots, A_{\epsilon, m}]$, where 
      $$ 
      A_{\epsilon, i}:=(P+\epsilon I_\cH)^{-1/2} X_i^* (P+\epsilon 
      I_\cH)^{-1/2}, 
      \quad 
      i=1,\ldots, m, 
      $$ 
      one can see that \eqref{pos1} holds if and only if $\|A_\epsilon\|\leq 1$ 
      for any $\epsilon >0$. 
      Since $V_1,\ldots, V_m$ are isometries with orthogonal ranges, we have $V_i^* V_j=\delta_{ij} I_\cK$ for $i,j=1,\ldots, m$, and
      \begin{equation*} 
      \begin{split} 
      \left\|\sum_{i=1}^m V_i \otimes A^*_{\epsilon, i}\right\|^2 &= 
      \left\| \sum_{i, j=1}^m V_i^* V_j \otimes A_{\epsilon, i}A^*_{\epsilon, j} 
      \right\|\\ 
      &=\left\|I_\cK\otimes \sum_{i=1}^m A_{\epsilon, i}A^*_{\epsilon, 
      i}\right\|= 
      \|A_\epsilon\|^2. 
      \end{split} 
      \end{equation*} 
      Therefore, $\sum\limits_{i=1}^m V_i \otimes A^*_{\epsilon, i}$ is a contraction 
      for 
      any $\epsilon >0$ if and only if $\|A_\epsilon\|\leq 1$ for any $\epsilon>0$. Applying again Lemma \ref{standard} to the   
      operator 
$\sum\limits_{i=1}^m V_i \otimes A^*_{\epsilon, i}$, 
      we obtain  $N(I_\cH, A^*_{\epsilon, i})\geq 0$
for any $\epsilon>0$, which is equivalent to 
      $N(P+\epsilon I_\cH, X_i)\geq 0$ for any $\epsilon>0$. Taking $\epsilon\to 
      0$, 
       we obtain 
      inequality \eqref{pos2}. The converse, follows exactly the same lines. 
      The proof is complete. 
      \end{proof}

       We need a few more definitions. 
       Given $\alpha, \beta\in \FF_n^+$, we say that $\alpha>\beta$ if 
        $\alpha=\beta \omega$ 
        for some  $\omega\in \FF_n^+\backslash\{g_0\}$. We denote 
        $\omega:=\alpha\backslash \beta$.  
        A kernel $K:\FF_n^+\times \FF_n^+\to B(\cH)$ is called multi-Toeplitz if 
        $K(g_0,g_0)=I_\cH$ and 
        $$ 
        K(\alpha, \beta)= 
        \begin{cases} 
          K(\alpha\backslash \beta, g_0) &\text{ if } \alpha>\beta\\ 
          I_\cH &\text{ if } \alpha=\beta\\ 
          K( g_0, \beta\backslash \alpha) &\text{ if } \alpha<\beta\\ 
          0\quad &\text{ otherwise}. 
         \end{cases}$$ 
        It is said to be positive definite provided that 
        $$ 
        \sum_{\alpha,\beta\in \FF_n^+} \left<K(\alpha,\beta) h(\beta), 
        h(\alpha)\right>\geq 0 
        $$ 
        for all finitely supported functions $h$ from $\FF_n^+$ into $\cH$.

 Using the noncommutative Poisson transforms  of \cite{Po-poisson}, one can easily  prove the following.

\begin{lemma}\label{pois} 
 Let 
       $$F(X_1,\ldots, X_n):=\sum\limits_{k=0}^\infty \sum\limits_{|\alpha|=k} 
      X_\alpha 
      \otimes A_{(\alpha)}, \quad A_{(\alpha)}\in B(\cH) 
      $$ be a  holomorphic function on $[B(\cX)^n]_{<1}$ with coefficients in $B(\cH)$.
Then the following statements are equivalent:
\begin{enumerate}
\item[(i)]$\text{\rm Re}\, F(X_1,\ldots,X_n)\leq I$ for any 
$(X_1,\ldots,X_n)\in [B(\cX)^n]_{<1}$ and any Hilbert space $\cX$;
\item[(ii)]$\text{\rm Re}\, F(rS_1,\ldots,r S_n)\leq I$  for any  $0\leq r<1$.
\end{enumerate}
\end{lemma}

 In what follows we denote by 
$[B_{(\alpha)}:\ |\alpha|=k]$ the row matrix with entries 
 $B_{(\alpha)}\in B(\cH)$, where $\alpha\in \FF_n^+$ and $|\alpha|=k$. 

      \begin{theorem}\label{noncom} 
       Let 
       $$F(X_1,\ldots, X_n):=\sum\limits_{k=0}^\infty \sum\limits_{|\alpha|=k} 
      X_\alpha 
      \otimes A_{(\alpha)}, \quad A_{(\alpha)}\in B(\cH), 
      $$ be a  holomorphic function on $[B(\cX)^n]_{<1}$ with coefficients in $B(\cH)$
such that  $F(0)\geq 0$   and       %
      $$  
      \text{\rm Re}\, F(X_1,\ldots,X_n)\leq I\quad  \text{ for any }  \ (X_1,\ldots, X_n)\in [B(\cX)^n]_{<1}. 
      $$
       Then 
      \begin{enumerate} 
      \item[(i)] the operator matrix 
      $$ 
      P_k:=\left( 
      \begin{matrix} 
      2(I_\cH -A_{(0)}) & [A^*_{(\alpha)}:\ |\alpha|=k]\\ 
      \left[\begin{matrix} A_{(\alpha)}\\ :\\|\alpha|=k \end{matrix}\right] 
      & \left[ \begin{matrix}2(I_\cH -A_{(0)})& \cdots & 0\\ 
      \vdots& \ddots & \vdots\\ 
      0& \cdots &  2(I_\cH -A_{(0)})\end{matrix}\right] 
      \end{matrix} 
      \right) 
      $$ 
      is positive for any $k=1,2,\ldots; $
      \item[(ii)] if $\{Y_{(\alpha)}\}_{\alpha\in \FF_n^+}$ is a sequence of 
      operators 
      in $B(\cK)$  with 
      $$ 
      \|Y_{(0)}\|\leq 1 \ \text{ and }\ 
      \sum_{k=1}^\infty \left\| \sum_{|\alpha|\\=k} 
      Y_{(\alpha)}^* Y_{(\alpha)}\right\|^{1/2}\leq \frac {1} {2}, 
      $$ 
      then 
      $$ 
      \left\|\sum_{k=0}^\infty \sum_{|\alpha|=k} S_\alpha \otimes  A_{(\alpha)} 
      \otimes Y_{(\alpha)}\right\|\leq 1, 
      $$ 
      where the series converges in the norm topology of 
      $B(F^2(H_n)\otimes \cH\otimes \cK)$. 
      \end{enumerate} 
      \end{theorem} 
      \begin{proof} Let $R_1,\ldots, R_n$ be the right creation operators 
       acting on the Fock 
       space $F^2(H_n)$. We recall that $R_i=U^* S_iU$, \ $i=1,\ldots, n$, 
       where $U$ is the flipping operator. Since $F(0)\geq 0$ and 
       $\text{\rm Re}\, F(rS_1,\ldots,r S_n)\leq I$  for $0\leq r<1$, we deduce 
       that $A_{(0)}\geq 0$ and 
       $$ 
       \sum_{k=1}^\infty \sum_{|\alpha|=k} r^{|\alpha|} R_\alpha^*\otimes 
      A_{(\alpha)}^* 
       + I \otimes 2A_{(0)} + 
       \sum_{k=1}^\infty \sum_{|\alpha|=k} r^{|\alpha|} R_\alpha\otimes 
      A_{(\alpha)} 
       \leq 2I 
       $$ 
       for  any  $0\leq r<1$, where the series are norm convergent. 
      Hence, we infer that 
      \begin{equation} \label{toeplitz} 
      \sum_{k=1}^\infty \sum_{|\alpha|=k}r^{|\alpha|} R_\alpha^*\otimes 
      C_{(\alpha)}^* 
       + I\otimes C_{(0)} + 
       \sum_{k=1}^\infty \sum_{|\alpha|=k} r^{|\alpha|}R_\alpha\otimes 
      C_{(\alpha)}\geq 0, 
      \end{equation} 
      where 
      \begin{equation}\label{condit} 
      C_{(0)}:= 2(I_\cH-A_{(0)}) \ \text{ and }\ C_{(\alpha)}:= -A_{(\alpha)} \ 
      \text{ if } \ \alpha \in \FF_n^+\backslash\{g_0\}. 
      \end{equation} 
      For each $r\in [0,1)$, define the multi-Toeplitz  kernel $K_{F,r}:\FF_n^+\times \FF_n^+\to 
      B(\cH)$ 
      by 
      \begin{equation}\label{ke} 
         K_{F,r}(\alpha, \beta):= 
         \begin{cases} 
         r^{|\beta\backslash \alpha|}C^*_{(\widetilde{\beta\backslash \alpha})} 
      &\text{ if } \beta>\alpha\\ 
         C_{(0)}  &\text{ if } \alpha=\beta\\ 
         r^{|\alpha\backslash \beta|} C_{(\widetilde{\alpha\backslash \beta})}  
      &\text{ if } \alpha>\beta\\ 
          0\quad &\text{ otherwise}, 
         \end{cases} 
         \end{equation} 
where $\tilde{\gamma}$ is the reverse of $\gamma\in \FF_n^+$.
       Note that if $\{h_\beta\}_{|\beta|\leq q}\subset \cH$, then 
          \begin{equation*} 
          \begin{split} 
          \left< \left(\sum_{k=0}^\infty \sum_{|\alpha|=k} 
          r^{|\alpha|} R_\alpha\otimes C_{(\alpha)}\right) 
           \right.& \left.\left(\sum_{|\beta|\leq q} e_\beta\otimes h_\beta 
      \right), 
         \sum_{|\gamma|\leq q} e_\gamma\otimes h_\gamma\right>\\ 
         &= 
           \sum_{k=0}^\infty \sum_{|\alpha|=k}\left<\sum_{|\beta|\leq q} 
          r^{|\alpha|} R_\alpha e_\beta\otimes C_{(\alpha)}h_\beta, 
          \sum_{|\gamma|\leq q} e_\gamma\otimes h_\gamma\right> \\ 
          &=\sum_{\alpha\in \FF_n^+}\sum_{|\beta|, |\gamma|\leq q} 
           r^{|\alpha|}\left< e_{\beta \tilde{\alpha}}, e_\gamma\right> 
          \left<C_{(\alpha)} h_\beta, h_\gamma\right>\\ 
          &= 
          \sum_{ \gamma\geq\beta; ~|\beta|, |\gamma|\leq q} 
          r^{|\gamma\backslash \beta|} \left<C_{(\widetilde{\gamma\backslash 
      \beta})} 
      h_\beta, h_\gamma\right>\\ 
           &= 
           \sum_{\gamma\geq\beta; ~|\beta|, |\gamma|\leq q}\left<K_{F,r} 
            (\gamma, \beta)h_\beta, h_\gamma\right>. 
          \end{split} 
          \end{equation*} 
         Hence, taking into account that $K_{F,r} 
            (\gamma, \beta)=K_{F,r}^* 
            ( \beta, \gamma)$ and   inequality \eqref{toeplitz}, we deduce that 
             $\left[K_{F,r} 
      (\alpha, \beta)\right]_{|\alpha|,|\beta|\leq q}\geq 0$ for 
      any 
              $r\in [0,1)$. 
          Taking $r\to 1$, we obtain 
          $\left[K_{F,1} 
      (\alpha, \beta)\right]_{|\alpha|,|\beta|\leq q}\geq 0$. 
          According to Theorem 3.1 of \cite{Po-posi} and using \eqref{ke}, we deduce that there is a completely 
      positive 
      linear map 
          $\mu:C^*(S_1,\ldots, S_n)\to B(\cH)$ such that 
          $$ 
          \mu(S_\alpha)=K_{F,1}(g_0, \alpha)=C^*_{(\widetilde \alpha)},\quad 
           \alpha\in \FF_n^+. 
          $$ 
         Using Stinespring's  representation theorem (see \cite{St}), we find a Hilbert space 
      $\cG\supseteq \cH$, 
         a $*$-representation $\pi: C^*(S_1,\ldots, S_n)\to B(\cG)$, and 
         a bounded operator $X:\cH\to \cG$ such that 
         $$ 
         \mu(f)=X^* \pi(f) X,\quad  f\in C^*(S_1,\ldots, S_n). 
         $$ 
         Denote $V_i:= \pi(S_i)$, \ $i=1,\ldots, n$, and notice that 
         \begin{equation}\label{eq1} 
         X^* V_\alpha X=\mu(S_\alpha)=C^*_{(\widetilde \alpha)}\quad 
         \text{ if  } \alpha\in \FF_n^+\backslash\{g_0\}, 
         \end{equation} 
         and $X^*X =\mu(I)=C_{(0)}$. 
        On the other hand, if $T:\cM\to \cG$ is a contraction, 
        then 
        $\left(\begin{matrix} 
        0&T\\T^*&0 
        \end{matrix} \right) 
        $ 
        is a selfadjoint contraction acting on the Hilbert space $\cG\oplus 
      \cM$. 
        Hence,  
        $\left(\begin{matrix} 
        I_\cG &-T\\-T^*& I_\cM 
        \end{matrix} \right) 
        $ 
        is a positive operator. Since  $V_1,\ldots, V_n\in B(\cG)$ are 
        isometries with orthogonal ranges, the operator $T:=[V_\alpha:\ 
      |\alpha|=k]$ 
         is a row isometry acting from  $\oplus_{i=1}^{n^k} \cG$ to $\cG$. 
         Using the above-mentioned result, we deduce that 
         $$ 
         \left( 
      \begin{matrix} 
      XX^* & [-XV_{\alpha}X^*:\ |\alpha|=k]\\ 
      \left[\begin{matrix} -XV^*_{\alpha}X^*\\ :\\|\alpha|=k \end{matrix}\right] 
      & \left[ \begin{matrix} XX^*& \cdots & 0\\ 
      \vdots& \ddots & \vdots\\ 
      0& \cdots &  XX^*\end{matrix}\right] 
      \end{matrix} 
      \right)\geq 0. 
         $$ 
          Hence and using relations \eqref{condit} and  \eqref{eq1}, we infer 
      that 
          $P_k\geq 0$, 
          which proves part (i) of the theorem. 
          
          To prove  (ii), denote 
          $d_k:=\left\| [ Y_{(\alpha)}^*:\ |\alpha|=k]\right\|$ \ 
          for $k=1,2,\ldots$. Lemma \ref{standard} 
          shows that 
          $$ 
          Q_k:= 
          \left( 
      \begin{matrix} 
      d_kI_\cK & [Y_{(\alpha)}^*:\ |\alpha|=k]\\ 
      \left[\begin{matrix} Y_{(\alpha)}\\ :\\|\alpha|=k \end{matrix}\right] 
      & \left[ \begin{matrix} d_kI_\cK& \cdots & 0\\ 
      \vdots& \ddots & \vdots\\ 
      0& \cdots &  d_kI_\cK\end{matrix}\right] 
      \end{matrix} 
      \right)\geq 0 
         $$
for any $k=1,2,\ldots$. 
         Since $P_k\geq 0$ and $Q_k\geq 0$, we have $P_k\otimes Q_k\geq 0$. 
          Compressing the operator matrix $P_k\otimes Q_k$ to appropriate 
      entries, we 
          deduce that the operator 
          \begin{equation}\label{bigmat} 
          \left( 
      \begin{matrix} 
      2d_k(I_\cH -A_{(0)})\otimes I_\cK & 
      [A^*_{(\alpha)}\otimes Y_{(\alpha)}^*:\ |\alpha|=k]\\ 
      \left[\begin{matrix} A_{(\alpha)}\otimes Y_{(\alpha)}\\ :\\|\alpha|=k 
       \end{matrix}\right] 
      & \left[ \begin{matrix}2d_k(I_\cH -A_{(0)})\otimes I_\cK& \cdots & 0\\ 
      \vdots& \ddots & \vdots\\ 
      0& \cdots &  2d_k(I_\cH -A_{(0)})\otimes I_\cK\end{matrix}\right] 
      \end{matrix} 
      \right) 
          \end{equation} 
       is positive for any $k=1,2,\ldots$.  
         Applying Lemma \ref{pos1} to the operator  matrix \eqref{bigmat}, 
         we deduce that  
         \begin{equation}\label{mat1} 
         \left(\begin{matrix} 
         I_{F^2(H_n)}\otimes 2d_k(I_\cH -A_{(0)})\otimes I_\cK & 
         \sum\limits_{|\alpha|=k} S_\alpha\otimes A_{(\alpha)} \otimes 
      Y_{(\alpha)}\\ 
         \sum\limits_{|\alpha|=k} S^*_\alpha\otimes A^*_{(\alpha)} \otimes 
      Y^*_{(\alpha)} & 
         I_{F^2(H_n)}\otimes 2d_k(I_\cH -A_{(0)})\otimes I_\cK 
         \end{matrix} \right) 
         \end{equation} 
        is positive for any $k=1,2,\ldots$, where $S_1,\ldots, S_n$ are the left 
        creation operators 
        on the Fock space $F^2(H_n)$. 
Hence, using inequality $I-A_0\leq I$ and Lemma \ref{standard}, we deduce that
$$\left\|\sum_{|\alpha|=k} S_\alpha \otimes A_{(\alpha)}\otimes Y_{(\alpha)}\right\|\leq 2d_k
$$
for $k=1,2,\ldots$. Therefore, the series 
$\sum\limits_{k=1}\sum\limits_{|\alpha|=k} S_\alpha \otimes A_{(\alpha)}\otimes Y_{(\alpha)}$ is convergent in norm.
        Since $A_{(0)}\geq 0$ and 
       $\left(\begin{matrix} 
        I_\cK &Y_{(0)}\\Y_{(0)}^*& I_\cK 
        \end{matrix} \right)\geq 0, 
        $ 
        we have 
        \begin{equation}\label{Azero} 
        \left(\begin{matrix} 
        I_{F^2(H_n)}\otimes A_{(0)}\otimes I_\cK & I_{F^2(H_n)}\otimes 
      A_{(0)}\otimes 
      Y_{(0)}\\ 
        I_{F^2(H_n)}\otimes A_{(0)}^*\otimes Y_{(0)}^*& I_{F^2(H_n)}\otimes 
      Y_{(0)}\otimes I_\cK 
        \end{matrix} \right)\geq 0. 
        \end{equation} 

        Taking the sum of the operator matrices given by \eqref{Azero} for $k=1,2,\ldots$,  and \eqref{mat1}, 
        and taking   into account that $\sum\limits_{k=1}^\infty d_k\leq \frac 
      {1} 
      {2}$, 
         we deduce that 
        $$ 
        \left(\begin{matrix} 
         I_{F^2(H_n)\otimes \cH\otimes \cK}  & 
        \sum\limits_{k=0}^\infty \sum\limits_{|\alpha|=k} S_\alpha\otimes 
      A_{(\alpha)} 
      \otimes Y_{(\alpha)}\\ 
        \sum\limits_{k=0}^\infty \sum\limits_{|\alpha|=k} S^*_\alpha\otimes 
      A^*_{(\alpha)} \otimes Y^*_{(\alpha)} & 
          I_{F^2(H_n)\otimes \cH\otimes \cK} 
         \end{matrix} \right) 
         $$ 
        is a  positive operator acting on the Hilbert space 
        $F^2(H_n)\otimes \cH\otimes \cK$. Now, Lemma \ref{standard} implies 
        $$ 
        \left\|\sum_{k=0}^\infty \sum_{|\alpha|=k} S_\alpha\otimes A_{(\alpha)} 
        \otimes Y_{(\alpha)}\right\|\leq 1. 
        $$ 
       The proof is complete.  
      \end{proof} 

The next two theorems provide Wiener and Bohr type inequalities for holomorphic functions on $[B(\cX)^n]_{<1}$ with coefficients in $B(\cH)$.

      \begin{theorem} \label{bohr-gen} 
      Let 
       $$F(X_1,\ldots, X_n):=\sum\limits_{k=0}^\infty \sum\limits_{|\alpha|=k} 
      X_\alpha 
      \otimes A_{(\alpha)}, \quad A_{(\alpha)}\in B(\cH), 
      $$
 be  a  holomorphic function  on  $[B(\cX)^n]_{<1}$    such that  
 $F(0)\geq 0$   and       %
      $$  
      \text{\rm Re}\, F(X_1,\ldots,X_n)\leq I\quad  \text{ for any }  \ (X_1,\ldots, X_n)\in [B(\cX)^n]_{<1}. 
      $$
       Then 
      \begin{enumerate} 
      \item[(i)] if $Y_1,\ldots, Y_n\in B(\cK)$ are bounded operators such that 
      $\|[Y_1^*,\ldots, Y_n^*]\|\leq \frac {1} {3}$, then 
      $$\|F(S_1\otimes Y_1,\ldots, S_n\otimes Y_n)\|\leq 1; 
      $$ 
      \item[(ii)] 
      $\sum\limits_{|\alpha|=k} A_{(\alpha)}^* A_{(\alpha)}\leq 4(I-A_{(0)})$ 
      for any $k=1,2,\ldots;$ 
      \item[(iii)] 
      $\left\|\sum\limits_{|\alpha|=k} A_{(\alpha)}^* A_{(\alpha)}\right\|^{1/2} 
      \leq 2\|I-A_{(0)}\|$ with coefficients in $B(\cH)$
      for any $k=1,2,\ldots;$ 
      \item[(iv)] for $0\leq r<1$, we have 
      $$ 
      \sum_{k=0}^\infty r^k\left( 
      \sum\limits_{|\alpha|=k} A_{(\alpha)}^* A_{(\alpha)}\right)^{1/2}\leq 
      M(r)I, 
      $$ 
      where 
      \begin{equation}\label{mr} 
      M(r):= 
      \begin{cases} 
           1+\frac{r^2} {(1-r)^2} &\text{ if }  0\leq r\leq \frac {1} {2}\\ 
         \frac{2r} {1-r}   &\text{ if } \frac {1} {2}<r<1 ; 
         \end{cases} 
         \end{equation} 
      \item[(v)] if    $\|A_{(0)}\|<1$ and $0\leq r<1$, then 
      $$ 
      A_{(0)}+\sum_{k=1}^\infty \sum_{|\alpha|=k}  r^k A_{(\alpha)}^* 
      (I-A_{(0)})^{-1} A_{(\alpha)}\leq K(r)I, 
      $$ 
      where 
      \begin{equation}\label{kr} 
      K(r):= 
      \begin{cases} 
           1  &\text{ if }  0\leq r\leq \frac {1} {5}\\ 
         \frac{4r} {1-r}   &\text{ if } \frac {1} {5}<r<1 ; 
         \end{cases} 
         \end{equation} 
      \item[(vi)] for $0\leq r\leq \frac {1} {3}$, we have 
      $$ 
      \sum_{k=0}^\infty r^k 
      \left\|\sum\limits_{|\alpha|=k} A_{(\alpha)}^* A_{(\alpha)}\right\|^{1/2} 
      \leq \|A_{(0)}\|+\|I-A_{(0)}\|. 
      $$ 
      \end{enumerate} 
      \end{theorem} 
      \begin{proof} 
      Define $Y_{(\alpha)}:= Y_\alpha$, $\alpha\in \FF_n^+$, and notice that 
      $$ 
      \sum_{k=1}^\infty\left\| \sum_{|\alpha|=k} Y_{(\alpha)}^* Y_{(\alpha)} 
      \right\|^{1/2} 
      \leq \sum_{k=1}^\infty\left\| \sum_{i=1}^n Y_i^* Y_i\right\|^{k/2} 
      =\sum_{k=1}^\infty \frac {1} {3^k}=\frac {1} {2} 
      $$ 
      for $k=1,2,\ldots$. 
      Applying Theorem \ref{noncom} part (ii), we obtain 
      $$ 
      \left\|\sum_{k=0} \sum_{|\alpha|=k} S_\alpha\otimes Y_\alpha\otimes 
      A_{(\alpha)}\right\|\leq 1, 
      $$ 
       which implies (i). 
       To prove (ii) and (iii), note that relations \eqref{eq1} and 
      \eqref{condit} 
      imply 
       \begin{equation*} 
       \begin{split} 
       \sum_{|\alpha|=k}A_{\widetilde \alpha}^*A_{\widetilde \alpha}&= 
       \sum_{|\alpha|=k} X^* V_\alpha XX^* V_\alpha^* X\\ 
       &\leq \|X\|^2X^*\left(\sum_{|\alpha|=k}V_\alpha V_\alpha^*\right) X\leq 
      \|X\|^2 
      X^*X\\ 
       &=4\|I-A_{(0)}\| (I-A_{(0)}). 
       \end{split} 
       \end{equation*} 
       Therefore, we have 
       \begin{equation}\label{in1} 
        \sum_{|\alpha|=k}A_{ \alpha}^*A_{\alpha} 
        \leq 
        4\|I-A_{(0)}\| (I-A_{(0)}) 
       \end{equation} 
for $k=1,2,\ldots$.
       According to Theorem \ref{noncom}, we have $A_{(0)}\geq 0$ and $P_k\geq 
      0$. 
      This implies 
       $I_\cH-A_{(0)}\geq 0$ and $\|I_\cH-A_{(0)}\|\leq 1$. 
       One can easily see that inequality \eqref{in1} implies (ii) and (iii). 
       Since (ii) implies 
       $$\left( \sum_{|\alpha|=k}A_{\alpha}^*A_{ \alpha}\right)^{1/2} 
       \leq 2(I-A_{(0)})^{1/2}, \quad k=1,2,\ldots,
       $$ 
       we deduce that 
       \begin{equation*} 
       \begin{split} 
       \sum_{k=0}^\infty r^k\left( 
      \sum\limits_{|\alpha|=k} A_{(\alpha)}^* A_{(\alpha)}\right)^{1/2}&\leq 
      A_{(0)}+\frac{2r} {1-r}(I-A_{(0)})^{1/2}\\ 
      &\leq \sup_{0\leq x\leq 1} \left\{x+\frac{2r} {1-r}\sqrt{1-x}\right\} I= 
      M(r)I, 
       \end{split} 
       \end{equation*} 
       where $M(r)$ is given by \eqref{mr}. 

       Now, assume that $\|A_{(0)}\|<1$. According to Lemma \ref{posi}, 
       the positivity of the operator $P_k$ is equivalent to the positivity of 
      the 
      operator 
       $$ 
       Q_k:=\left(\begin{matrix} 
         I_{F^2(H_n)}\otimes 2(I_\cH -A_{(0)})  & 
         \sum\limits_{|\alpha|=k} S_\alpha\otimes A_{(\alpha)}  \\ 
         \sum\limits_{|\alpha|=k} S^*_\alpha\otimes A^*_{(\alpha)}  & 
         I_{F^2(H_n)}\otimes 2(I_\cH -A_{(0)}) 
         \end{matrix} \right) 
         $$ 
          for any $ k=1,2,\ldots$. 
         Using Lemma \ref{standard}, we deduce that $Q_k\geq 0$ if and only if 
         $$ 
         \left( \sum\limits_{|\alpha|=k} S^*_\alpha\otimes A^*_{(\alpha)}\right) 
         [I\otimes (I-A_{(0)})^{-1}] 
         \left( \sum\limits_{|\alpha|=k} S_\alpha\otimes A_{(\alpha)}\right) 
         \leq I\otimes4 (I-A_{(0)}) 
         $$ 
         for any $ k=1,2,\ldots$. 
         Hence, and taking into account that $S_i^* S_j=\delta_{ij} I$, \ 
      $i,j=1,\ldots, n$, 
         we obtain 
         $$\sum_{|\alpha|=k} A_{(\alpha)}^* (I-A_{(0)})^{-1} A_{(\alpha)}\leq 
      4(I-A_{(0)}), \quad k=1,2,\ldots.
         $$ 
         For $0\leq r<1$, the latter inequality implies 
         \begin{equation*} 
         \begin{split} 
         A_{(0)}+\sum_{k=1}^\infty \sum_{|\alpha|=k}  r^k A_{(\alpha)}^* 
      (I-A_{(0)})^{-1} A_{(\alpha)}&\leq 
      A_{(0)} +\frac {4r} {1-r} (I-A_{(0)})\\ 
      &\leq 
      \sup_{0\leq x\leq 1}\left\{ x+\frac {4r} {1-r}(1-x)\right\}I 
      \leq K(r)I, 
         \end{split} 
         \end{equation*} 
         where $K(r)$ is given by \eqref{kr}. 
         
         To prove (vi), notice that (iii) implies 
         \begin{equation*}\begin{split} 
         \sum_{k=0}^\infty r^k 
      \left\|\sum\limits_{|\alpha|=k} A_{(\alpha)}^* A_{(\alpha)}\right\|^{1/2} 
      &\leq \|A_{(0)}\|+2\frac {r} {1-r} \|I-A_{(0)}\|\\ 
      &\leq 
      \|A_{(0)}\|+\|I-A_{(0)}\| 
         \end{split} 
         \end{equation*} 
       for $0\leq r\leq \frac {1} {3}$. 
       The proof is complete.  
      \end{proof}

Now we can prove the following operator-valued Bohr inequality for holomorphic functions on $[B(\cX)^n]_{<1}$.

\begin{theorem}
    \label{Bo-ge} 
       Let 
       $$F(X_1,\ldots, X_n):=\sum\limits_{k=0}^\infty \sum\limits_{|\alpha|=k} 
      X_\alpha 
      \otimes A_{(\alpha)}, \quad A_{(\alpha)}\in B(\cH), 
      $$
 be  a  holomorphic function  on  $[B(\cX)^n]_{<1}$ with coefficients in $B(\cH)$   such that  
 $F(0)\geq 0$   and       %
      $$  
      \text{\rm Re}\, F(X_1,\ldots,X_n)\leq I\quad  \text{ for any }  \ (X_1,\ldots, X_n)\in [B(\cX)^n]_{<1}. 
      $$ 
      Then 
      \begin{equation*} 
      \sum_{k=0}^\infty \left\| \sum_{|\alpha|=k} T_\alpha \otimes 
      A_{(\alpha)}\right\| 
      \leq \|A_{(0)}\|+\|I-A_{(0)}\| 
      \end{equation*} 
      for  any 
      $(T_1,\ldots, T_n)\in [B(\cK)^n]_{1/3}$ and any Hilbert space $\cK$. 
      \end{theorem} 
      \begin{proof} 
      Let $(T_1,\ldots, T_n)\in [B(\cK)^n]_{r}$ with $0\leq r\leq \frac {1} {3}$. 
      Notice that $[r^{-1}T_1,\ldots, r^{-1} T_n]$ is a row contraction and, 
      according 
      to 
       the noncommutative von Neumann inequality \cite{Po-von}, we have 
       \begin{equation*} 
       \begin{split} 
       \sum_{k=0}^\infty \left\| \sum_{|\alpha|=k} T_\alpha \otimes 
      A_{(\alpha)}\right\| 
       &\leq 
       \sum_{k=0}^\infty r^k\left\| \sum_{|\alpha|=k} S_\alpha \otimes 
        A_{(\alpha)}\right\| 
       \\ 
       &= 
       \sum_{k=0}^\infty r^k\left\| \sum_{|\alpha|=k} A_{(\alpha)}^* 
        A_{(\alpha)}\right\|^{1/2} 
       \end{split} 
       \end{equation*} 
       Using  now Theorem \ref{bohr-gen} part (vi), we can complete the proof. 
      \end{proof} 

 For each $\alpha\in \FF_n^+$,  define 
      $e_\alpha :=  e_{i_1}\otimes e_{i_2}\otimes \cdots \otimes e_{i_k}$  
       and $e_{g_0}= 1$. 
      It is  clear that 
      $\{e_\alpha:\alpha\in\FF_n^+\}$ is an orthonormal basis of $F^2(H_n)$. 
       We denote by $\cP$ the set of all polynomials in $F^2(H_n)$, 
       i.e., all the elements of the form 
       $$p=\sum_{|\alpha|\leq m } a_\alpha e_\alpha, 
       \quad a_\alpha\in \CC,  \ m=0,1,2,\ldots. 
$$

\begin{remark}
 Let 
      $$F(S_1,\ldots, S_n):=\sum\limits_{k=0}^\infty \sum\limits_{|\alpha|=k} 
      S_\alpha 
      \otimes A_{(\alpha)}, \quad A_{(\alpha)}\in B(\cH), 
      $$ be  a formal  power series  such that,  
for each \ $h\in \cH$,
$$
\limsup_{k\to\infty}
\left(\sum_{|\alpha|=k} \|A_\alpha h\|^2\right)^{1/2k}\leq 1.
$$ 
If 
            $F(0)\geq 0$   and      
     $$ \text{\rm Re}\, \left<F(rS_1,\ldots, rS_n)p,p\right> \leq \|p\|^2\quad
\text{ for any } \ p\in \cP\otimes \cH \ \text{ and  } \ 0\leq r<1,   $$
  then the conclusions of 
Theorem $\ref{noncom}$, Theorem $\ref{bohr-gen}$, and Theorem $\ref{Bo-ge}$ remain true.
\end{remark}

      We need to recall from  \cite{Po-charact}, 
      \cite{Po-multi}, \cite{Po-von},  \cite{Po-funct}, and  \cite{Po-analytic} 
       a few facts 
       concerning multi-analytic   operators on Fock 
      spaces. 
         We say that 
       a bounded linear 
        operator 
      $M$ acting from $F^2(H_n)\otimes \cK$ to $ F^2(H_n)\otimes \cK'$ is 
       multi-analytic 
      if 
      \begin{equation*} 
      M(S_i\otimes I_\cK)= (S_i\otimes I_{\cK'}) M\quad 
      \text{\rm for any }\ i=1,\dots, n. 
      \end{equation*} 
      Notice that $M$ is uniquely determined by   the ``coefficients'' 
        $\theta_{(\alpha)}\in B(\cK, \cK')$  given by 
       $$ 
      \left< \theta_{(\tilde\alpha)}k,k'\right>:=  \left< M(1\otimes k), e_\alpha 
      \otimes k'\right>,\quad 
      k\in \cK,\ k'\in \cK',\ \alpha\in \FF_n^+, 
      $$ 
      where $\widetilde\alpha$ is the reverse of $\alpha$, i.e., 
      $\widetilde\alpha= 
      g_{i_k}\cdots g_{i_1}$ if 
      $\alpha= g_{i_1}\cdots g_{i_k}$. 
We denote $\theta_{(0)}:=\theta_{(g_0)}$.
      Note that 
      $$ 
      \sum\limits_{\alpha \in \FF_n^+} \theta_{(\alpha)}^*   
      \theta_{(\alpha)}\leq 
      \|M\|^2 I_\cK. 
      $$ 
       We can associate with $M$ a unique formal Fourier expansion 
      \begin{equation*}       M\sim \sum_{\alpha \in \FF_n^+} R_\alpha \otimes \theta_{(\alpha)}, \end{equation*} 
      where $R_i:=U^* S_i U$, \ $i=1,\ldots, n$, are the right creation 
      operators 
      on $F^2(H_n)$ and 
      $U$ is the (flipping) unitary operator on $F^2(H_n)$ mapping 
        $e_{i_1}\otimes e_{i_2}\otimes\cdots\otimes e_{i_k}$ into 
       $e_{i_k}\otimes\cdots\otimes e_{i_2}\otimes e_{i_1} $. 
       Since  the operator $M$ acts like its Fourier representation on ``polynomials'', 
  we will identify them for simplicity. 
      Based on the noncommutative von Neumann 
         inequality (\cite{Po-von}, \cite{Po-funct}),         we proved that 
        $$M=\text{\rm SOT}-\lim_{r\to 1}\sum_{k=0}^\infty 
      \sum_{|\alpha|=k} 
         r^{|\alpha|} R_\alpha\otimes \theta_{(\alpha)}, 
         $$ 
         where, for each $r\in [0,1)$, the series converges in the uniform norm. 
      Moreover, the set of  all multi-analytic operators in 
      $B(F^2(H_n)\otimes \cK, 
      F^2(H_n)\otimes \cK')$  coincides  with   
      $R_n^\infty\bar \otimes B(\cK,\cK')$, 
      the WOT closed algebra generated by the spatial tensor product, where 
      $R_n^\infty=U^* F_n^\infty U$. 

Now we can deduce the following operator-valued generalization of Bohr's type inequality for the noncommuative analytic Toeplitz algebra $F_n^\infty$. 

      \begin{corollary}\label{Bo-co} 
     If $F(S_1,\ldots, S_n):=\sum\limits_{k=0}^\infty \sum\limits_{|\alpha|=k} S_\alpha 
      \otimes A_{(\alpha)}$ is  in $ F_n^\infty\bar \otimes B(\cH)$, 
      then 
      $$ 
      \sum_{k=0}^\infty \left\| \sum_{|\alpha|=k} T_\alpha \otimes 
      A_{(\alpha)}\right\| 
      \leq \|F(S_1,\ldots, S_n)\|(\|A_{(0)}\|+\|I-A_{(0)}\|) 
      $$ 
      for  any 
      $(T_1,\ldots, T_n)\in [B(\cK)^n]_{1/3}$. 
Moreover, if $F(0)\geq 0$ and $\text{\rm Re}\, F(S_1,\ldots, S_n)\leq I$, then 
$$ 
      \sum_{k=0}^\infty \left\| \sum_{|\alpha|=k} T_\alpha \otimes 
      A_{(\alpha)}\right\| 
      \leq (\|A_{(0)}\|+\|I-A_{(0)}\|) 
      $$ 
      for  any 
      $(T_1,\ldots, T_n)\in [B(\cH)^n]_{1/3}$. 

      \end{corollary} 
      \begin{proof} Without loss of generality, we can assume that 
      $ \|F(S_1,\ldots, S_n)\|=1$.   Now, one can use the
       noncommutative Poisson 
      transforms of \cite{Po-poisson} to show that the hypotheses of Theorem \ref{Bo-ge} are 
      satisfied. 
       Applying the latter theorem, the result follows. 
      \end{proof} 

We remark that, in the particular case when $\cH=\CC$, the first part of Corollary \ref{Bo-co} was obtained in \cite{PPoS}. Its second part provides an operator-valued extension of Corollary \ref{Boh3} as well as a new proof.

We can improve some of the inequalities of Theorem \ref{bohr-gen} under the more restrictive  conditions that 
$F(S_1,\ldots, S_n)\in F_n^\infty\bar \otimes B(\cH)$ and 
$\|F(S_1,\ldots, S_n)\|\leq 1$. 

      \begin{theorem}\label{oper-gen} 
      If 
      $$ 
      F(S_1,\ldots, S_n):=\sum\limits_{k=0}^\infty \sum\limits_{|\alpha|=k} 
      S_\alpha 
      \otimes A_{(\alpha)}\in F_n^\infty\bar \otimes B(\cH) 
      $$ 
       and\  
      $\|F(S_1,\ldots, S_n)\|\leq 1$, then 
      \begin{enumerate} 
      \item[(i)] 
       $ 
      \left(\sum\limits_{|\alpha|=k} A_{(\alpha)}^*A_{(\alpha)}\right)^{1/2} 
      \leq (I-A_{(0)}^* A_{(0)})^{1/2}\quad 
      \text{ for any }\  k=1,2,\ldots; 
      $ 
      \item[(ii)] for $0\leq r<1$, we have 
      $$ 
      \sum_{k=0}^\infty r^k\left( 
      \sum\limits_{|\alpha|=k} A_{(\alpha)}^* A_{(\alpha)}\right)^{1/2}\leq 
      \left( 1+\frac{r^2} {(1-r)^2}\right)^{1/2}I; 
      $$ 
      \item[(iii)] if    $\|A_{(0)}\|<1$ and $k=1,2,\ldots$, then 
      $$ 
      \sum_{|\alpha|=k}  A_{(\alpha)}^* 
      (I-A_{(0)} A^*_{(0)})^{-1} A_{(\alpha)}\leq I-A_{(0)}^*A_{(0)}; 
      $$ 
      \item[(iv)] 
      if    $\|A_{(0)}\|<1$, $A_{(0)}\geq 0$,  and $0\leq r\leq 1$, then 
      $$ 
      A_{(0)}^2+\sum_{k=1}^\infty \sum_{|\alpha|=k}  r^k A_{(\alpha)}^* 
      (I- A^2_{(0)})^{-1} A_{(\alpha)}\leq N(r) I, 
      $$ 
      where 
      \begin{equation}\label{nr} 
      N(r):= 
      \begin{cases} 
           1  &\text{ if }  0\leq r\leq \frac {1} {2}\\ 
         \frac{r} {1-r}   &\text{ if } \frac {1} {2}<r<1. 
         \end{cases} 
         \end{equation} 
      \end{enumerate} 
      \end{theorem} 
      \begin{proof} 
      Let $\cM$ be the subspace of $F^2(H_n)$ spanned by the vectors $1$ and 
      $e_\alpha $, where $\alpha\in \FF_n^+$ and  
      $|\alpha|=k$. The operator $P_{\cM\otimes \cH} F(S_1,\ldots, 
      S_n)|\cM\otimes 
      \cH$ 
      is a contraction and its operator matrix with respect to the decomposition 
      $\cM\otimes \cH=\cH\oplus \bigoplus\limits_{|\alpha|=k}(e_\alpha\otimes 
      \cH)$ is 
      $$ 
      Y_k= 
      \left(\begin{matrix} A_{(0)} & [\begin{matrix}0 &\cdots &0\end{matrix}]\\ 
      \left[ \begin{matrix} 
      A_{(\alpha)}\\:\\|\alpha|=k\end{matrix}\right] 
      & \left[ \begin{matrix} 
      A_{(0)}&\cdots & 0\\ \vdots&\ddots & \vdots\\ 
      0&\cdots & A_{(0)}\end{matrix}\right] 
      \end{matrix} \right). 
      $$ 
      Hence, 
      $ 
      \left(\begin{matrix} A_{(0)}  \\ 
      \left[ \begin{matrix} 
      A_{(\alpha)}\\:\\|\alpha|=k\end{matrix}\right] 
      \end{matrix} \right) 
      $ 
      is a  contraction, which implies 
      $$ 
      \sum_{|\alpha|=k} A_{(\alpha)}^*A_{(\alpha)} 
      \leq I-A_{(0)}^* A_{(0)} \quad \text{ for }\ k=1,2\ldots. 
      $$ 
      The latter inequality implies (i). Using (i), we infer that 
      \begin{equation*} 
      \begin{split} 
      \sum_{k=0}^\infty r^k\left( 
      \sum\limits_{|\alpha|=k} A_{(\alpha)}^* A_{(\alpha)}\right)^{1/2}&\leq 
      (A_{(0)}^* A_{(0)})^{1/2} +\frac {r} {1-r} (1-A_{(0)}^* A_{(0)})^{1/2}\\ 
      &\leq \sup_{0\leq x\leq 1} \left\{ x+\frac {r} {1-r}\sqrt{1-x^2}\right\} 
      I\\ 
      &= 
      \left( 1+\frac{r^2} {(1-r)^2}\right)^{1/2}I. 
      \end{split} 
      \end{equation*} 
      Now, we prove part (iii). Taking into account the  structure of $2\times 
      2$ 
      lower 
       triangular 
      contractions, the operator 
      $$\left[ 
      \begin{matrix} 
      A&0\\C&D 
      \end{matrix}\right]:\cH\oplus \cK\to \cH'\oplus \cK' 
      $$ 
      is a contraction if and  only if $A$ and $D$ are contractions and 
      $$ 
      C=(I-DD^*)^{1/2} \Gamma (I-A^*A)^{1/2}, 
      $$ 
      where $\Gamma:\cD_A\to \cD_{D^*}$ is a contraction,  and the subspaces 
      $\cD_A$ 
      and $\cD_{D^*}$ 
      are defined by 
      $\cD_A:=\overline{(I-A^*A)^{1/2}\cH}$ and 
      $\cD_{D^*}:=\overline{(I-DD^*)^{1/2}\cK'}$. 
      Hence, if $\|D\|<1$, then we have 
      $$ 
      C^*(I-DD^*)^{-1} C=(I-A^*A)^{1/2}\Gamma^* \Gamma (I-A^*A)^{1/2}\leq 
      I-A^*A. 
      $$ 
     Applying   this result to the matrix $Y_k$, where $\|A_{(0)}\|<1$, we obtain 
      $$ 
      [A_{(\alpha)}^*:\ |\alpha|=k] 
      \left[\begin{matrix} 
      (I-A_{(0)} A_{(0)}^*)^{-1}&\cdots& 0\\ 
      \vdots&\ddots&\vdots\\ 
      0&\cdots& (I-A_{(0)} A_{(0)}^*)^{-1} 
      \end{matrix} 
      \right] 
      \left[ 
      \begin{matrix} 
      A_{(\alpha)}\\ 
      :\\ 
      |\alpha|=k 
      \end{matrix} 
      \right]\leq I-A_{(0)}^* A_{(0)} 
      $$ 
for $k=1,2,\ldots$,
      which proves (iii). 
      Now, assume that $A_{(0)}\geq 0$ and $\|A_{(0)}\|<1$. Using (iii), we 
      deduce 
      that 
      \begin{equation*} 
      \begin{split} 
      A_{(0)}^2+\sum_{k=1}^\infty \sum_{|\alpha|=k}  r^k A_{(\alpha)}^* 
      (I- A^2_{(0)})^{-1} A_{(\alpha)} 
      &\leq 
      A_{(0)}^2+\frac {r} {1-r}(I-A_{(0)}^2)\\ 
      &\leq  \sup_{0\leq x\leq 1} \left\{ x+\frac {r} {1-r}(1-x)\right\} 
      I=N(r)I, 
      \end{split} 
      \end{equation*} 
      where $N(r)$ is given by \eqref{nr}. 
      The proof is complete. 
      \end{proof} 

 We   recall  from \cite{Po-unitary} the following  multivariable operator-valued generalization of
     the  inequalities
  of Fej\' er and Egerv\'ary-Sz\'azs, to the spatial tensor product $C^*(S_1,\ldots, S_n)\otimes B(\cH)$. 
   Let $m\geq 2$   and let  $\left\{A_{(\alpha)}\right\}_{|\alpha|\leq m-1}$  be 
a sequence of operators in $B(\cH)$ such that
  the  operator
 \begin{equation*}
 \sum_{1\leq k\leq m-1}  S_\alpha^*\otimes A_{(\alpha)}+  
 I\otimes A_{0}+
 \sum_{1\leq k\leq m-1}  S_\alpha \otimes A_{(\alpha)}^*
 \end{equation*}
is positive.  
 Then, 
\begin{equation}\label{joint-ine}
w_e\left(A_{(\alpha)}:\ |\alpha|=k\right)
\leq w\left(A_{(\alpha)}:\ |\alpha|=k\right)
\leq 
 \|A_{0}\|\cos\frac {\pi} {\left[\frac{m-1} {k}\right]+2}
\end{equation}
 for  $1\leq k\leq m-1$,
where  $[x]$ is the integer part of $x$,  $w\left(X_1,\ldots, X_n\right)$ is the joint numerical radius of the $n$-tuple $(X_1,\ldots, X_n)\in B(\cX)^n$, i.e.,
$$ w\left(X_1,\ldots, X_n\right):=\omega(S_1\otimes X_1^*+\cdots +S_n\otimes X_n^*),
$$
and
  $w_e\left(X_1,\ldots, X_n\right)$ is the euclidean  joint numerical radius of the $n$-tuple $\left(X_1,\ldots, X_n\right)$, i.e.,
$$w_e\left(X_1,\ldots, X_n\right):= \sup_{\|h\|=1}\left( \sum_{i=1}^n|\left< X_ih,h\right>|^2\right)^{1/2}.
 $$
We also recall  that both the joint numerical radius and the euclidean joint numerical radius  are norms equivalent to
the operator norm on $B(\cX)^n$. Moreover,
$$
\frac{1}{2} \|[X_1,\ldots, X_n]\|\leq w\left(X_1,\ldots, X_n\right)\leq \|[X_1,\ldots, X_n]\|
$$
and 
$$
\frac{1}{2\sqrt{n}} \|[X_1,\ldots, X_n]\|\leq w_e\left(X_1,\ldots, X_n\right)
\leq w\left(X_1,\ldots, X_n\right).
$$

In what follows we obtain an operator-valued Bohr type inequality when the norm of the coefficients is replaced by the joint numerical radius. The result is new even in the single variable case $n=1$.

\begin{theorem}\label{pol} Let  $$F(S_1,\ldots, S_n):= \sum_{k=0}^{m-1}\sum_{|\alpha|=k}S_\alpha\otimes A_{(\alpha)}, \quad 
A_{(\alpha)}\in B(\cH),$$ 
be a polynomial 
 such that \ $F(0)\geq 0$\  and  
$ \text{\rm Re}\,F(S_1,\ldots, S_n)\leq I$. Then 
$$
\sum_{k=0}^{m-1}r^k w\left(A_{(\alpha)}^*:\ |\alpha|=k\right)\leq  \|A_{(0)}\|+\|I-A_{(0)}\|
$$
for any  $r\in [0,t_m]$, where $t_m\in (0,1]$ is the solution of the equation 
\begin{equation*}
\sum_{k=1}^{m-1} t^k\cos\frac{\pi} {\left[\frac{m-1}{k}\right]+2}=\frac{1} {2},
\end{equation*}
where $[x]$ is the integer part of $x$. Moreover $\{t_m\}$ is a strictly decreasing sequence which converges to $\frac{1}{3}$.
\end{theorem}
\begin{proof}
The conditions $F(0)\geq 0$ and 
 $\text{\rm Re}\,F(S_1,\ldots, S_n)\leq I$ imply
$$
\sum_{1\leq |\alpha|\leq m-1} -  S_\alpha^*\otimes A_{(\alpha)}^* +I\otimes 2(I-A_0)+ \sum_{1\leq |\alpha|\leq m-1} -  S_\alpha\otimes A_{(\alpha)}\geq 0.
$$
According to the inequality 
\eqref{joint-ine}, we have
\begin{equation*}
  w\left(A^*_{(\alpha)}:\ |\alpha|=k\right)
\leq 
 2\|I-A_{0}\|\cos\frac {\pi} {\left[\frac{m-1} {k}\right]+2}
\end{equation*}
 for  $1\leq k\leq m-1$.
If $0\leq r\leq t_m$, then  
we have 
   \begin{equation*}\begin{split}
\sum_{k=0}^{m-1}r^k w\left(A_{(\alpha)}^*:\ |\alpha|=k\right)&\leq 
 \|A_0\|+2\|I-A_0\|\sum_{k=1}^{m-1} r^k\cos\frac{\pi} {\left[\frac{m-1}{k}\right]+2}\\
&\leq \|A_0\|+2\|I-A_0\|\sum_{k=1}^{m-1} t_m^k\cos\frac{\pi} {\left[\frac{m-1}{k}\right]+2}\\
&\leq  \|A_{(0)}\|+\|I-A_{(0)}\|.
  \end{split}
\end{equation*}
The last part of the theorem was proved in Theorem \ref{sharp}.
The proof is complete.
 \end{proof}

\bigskip

\section{Operator-valued Bohr inequalities for   harmonic functions} 

We say that $G$ is a selfadjoint harmonic function on $[B(\cX)^n]_{<1}$ with coefficients in $B(\cH)$ if there exists a (universal) holomorphic function 
$F:[B(\cX)^n]_{<1}\to B(\cX\otimes \cH)$ such that
$$
G(X_1,\ldots, X_n)=\text{Re}\, F(X_1,\ldots, X_n)\quad \text{ 
for any }\ (X_1,\ldots, X_n)\in [B(\cX)^n]_{<1}.
$$
      Let $\cA:=\{A_{(\alpha)}\}_{\alpha\in \FF_n^+}$ be a sequence of operators 
      in $B(\cH)$  such that $A_{(0)}=A_{(0)}^*$ and 
$$ 
       F(X_1,\ldots, X_n)=\frac{1}{2}I\otimes A_{(0)}+ \sum_{k=1}^\infty \sum_{|\alpha|=k} 
        X_\alpha \otimes A_{(\alpha)}
        $$
is a holomorphic function on $[B(\cX)^n]_{<1}$. Define 
$H_\cA:[B(\cX)^n]_{<1}\to B(\cH\otimes \cX)$  by setting
$H_\cA(X_1,\ldots, X_n)=\text{Re}\, F(X_1,\ldots, X_n)$
 for any $(X_1,\ldots, X_n)\in [B(\cX)^n]_{<1}$.

       We remark that, in the particular case when  
       $\cA:=\{a_{\alpha}\}_{\alpha\in \FF_n^+}\subset \CC$ are 
       the  coefficients of an element 
        in the noncommutative analytic Toeplitz algebra $F_n^\infty$ of the form  
       $$ 
       F(S_1,\ldots, S_n)=\frac{a_0} {2} I+ \sum_{k=1}^\infty \sum_{|\alpha|=k} 
        a_\alpha S_\alpha,\quad a_0=\bar a_0, 
        $$        
then    the operator 
       $$ 
       H_\cA(S_1,\ldots, S_n):=\frac{1}{2}[F(S_1,\ldots, S_n)^* + F(S_1,\ldots, S_n)] 
       $$ 
       can be seen as a noncommutative analogue of the boundary function of a  (real-valued) 
         bounded harmonic function in the unit disc, while  $H_\cA(X_1,\ldots, X_n)$ can be seen as the 
noncommutative Poisson transform of  $H_\cA(S_1,\ldots, S_n)$ at the point $(X_1,\ldots, X_n)\in[B(\cX)^n]_{<1}$.

We  remark that, using the noncommutative Poisson transforms of \cite{Po-poisson}, one can easily show that 
$$H_\cA(X_1,\ldots, X_n)\leq H_\cB(X_1,\ldots, X_n) \quad
\text{ for any } \ (X_1,\ldots, X_n)\in [B(\cX)^n]_{<1}
$$
 if and only if
$$H_\cA(rS_1,\ldots, rS_n)\leq H_\cB(rS_1,\ldots, rS_n)
\text{ for any } \ r\in [0,1).$$

The first result of this section provides  Wiener and Bohr type inequalities for the coefficients of two harmonic functions on $[B(\cX)^n]_{<1}$ satisfying the inequality
$H_\cA(X_1,\ldots, X_n) \leq H_\cB(X_1,\ldots, X_n)$.
      \begin{theorem}\label{f<g} 
      Let $\cA:=\{A_{(\alpha)}\}_{\alpha\in \FF_n^+}$  and 
      \ $\cB:=\{B_{(\alpha)}\}_{\alpha\in \FF_n^+}$  be  sequences of operators 
      in $B(\cH)$ such that $H_\cA$ and $H_\cB$ are noncommutative harmonic functions on $[B(\cX)^n]_{<1}$ and  
       \begin{equation}\label{hh} 
      H_\cA(X_1,\ldots, X_n) \leq H_\cB(X_1,\ldots, X_n)\quad \text{ for  any 
      }\ (X_1,\ldots, X_n)\in[B(\cX)^n]_{<1}. 
      \end{equation} 
       Then 
       \begin{equation}\label{BA} 
        \sum_{|\alpha|=k} (B_{(\alpha)}- A_{(\alpha)})^*(B_{(\alpha)}- 
      A_{(\alpha)}) 
         \leq \|B_{(0)}-A_{(0)}\| (B_{(0)}-A_{(0)}) 
       \end{equation} 
for $k=1,2,\ldots$,
       and 
       \begin{equation*}         \sum_{k=1}^\infty  r^k \left\| \sum_{|\alpha|=k} 
       A_{(\alpha)}^* A_{(\alpha)}\right\|^{1/2}\leq 
       \frac{1} {2} \|B_{(0)}-A_{(0)}\|  + 
      \sum_{k=1}^\infty r^k \left\| \sum_{|\alpha|=k} 
       B_{(\alpha)}^* B_{(\alpha)}\right\|^{1/2}  
       \end{equation*} 
for any $r\in \left[0,\frac{1}{3}\right]$.
      \end{theorem} 
        
        \begin{proof} Notice that the inequality \eqref{hh} implies 
        \begin{equation}\label{poz2} 
        \sum_{k=1}^\infty \sum_{|\alpha|=k} r^{|\alpha|} S_\alpha^*\otimes 
      C_{(\alpha)}^* 
       + I\otimes C_{(0)}+ 
       \sum_{k=1}^\infty \sum_{|\alpha|=k} r^{|\alpha|} S_\alpha\otimes 
      C_{(\alpha)}\geq 0,\quad r\in [0,1), 
        \end{equation} 
       where 
       $  C_{(\alpha)}:= B_{(\alpha)}- A_{(\alpha)}$ if $\alpha\in 
      \FF_n^+\backslash\{g_0\}$ 
       and $C_{(0)}:=B_{(0)}-A_{(0)}$.  First we consider the case when that $C_{(0)}\neq 
      0$. 
      According to 
       \eqref{poz2}, we have $C_{(0)}\geq 0$. As in the proof of Theorem 
      \ref{noncom} 
        part (i), 
       we deduce that 
       there is a completely positive linear map 
          $\mu:C^*(S_1,\ldots, S_n)\to B(\cH)$ such that 
          $$ 
          \mu(S_\alpha)=C^*_{(\widetilde \alpha)},\quad 
           \alpha\in \FF_n^+. 
          $$ 
         Using Stinespring's  representation theorem, we find a Hilbert space 
      $\cG\supseteq \cH$, 
         a $*$-representation $\pi: C^*(S_1,\ldots, S_n)\to B(\cG)$, and 
         a bounded operator $X:\cH\to \cG$ such that 
         $$ 
         \mu(f)=X^* \pi(f) X,\quad  f\in C^*(S_1,\ldots, S_n). 
         $$ 
         Denote $V_i:= \pi(S_i)$, \ $i=1,\ldots, n$, and notice that 
         \begin{equation}\label{eq11} 
         X^* V_\alpha X=\mu(S_\alpha)=C^*_{(\widetilde \alpha)}\quad 
         \text{ if  } \alpha\in \FF_n^+\backslash\{g_0\}, 
         \end{equation} 
         and $X^*X =\mu(I)=C_{(0)}$. 
       Hence, we infer that 
      \begin{equation*}\begin{split} 
      \sum_{|\alpha|=k} C_{\tilde{\alpha}}^*C_{\tilde{\alpha}}&= 
      \sum_{|\alpha|=k}X^* V_\alpha XX^* V_\alpha^* X\\ 
      &\leq \|X\|^2 X^*\left(\sum_{|\alpha|=k} V_\alpha V_\alpha^*\right) X\\ 
      &\leq \|X\|^2 X^*X=\|C_{(0)}\| C_{(0)}. 
      \end{split} 
      \end{equation*} 
      Therefore, the inequality \eqref{BA} holds. It is clear that the latter 
inequality implies 
 \begin{equation}\label{BA1} 
        \left\|\sum_{|\alpha|=k} (B_{(\alpha)}- A_{(\alpha)})^*(B_{(\alpha)}- 
      A_{(\alpha)})\right\|^{1/2} 
         \leq \|B_{(0)}-A_{(0)}\| 
       \end{equation} 
for $k=1,2,\ldots $.
Using this inequality, we deduce   that 
\begin{equation*}
\begin{split}
\sum_{k=1}^\infty r^k\left\|[A_{(\alpha)}:\ |\alpha|=k]\right\|&-\sum_{k=1}^\infty r^k\left\|[B_{(\alpha)}:\ |\alpha|=k]\right\|\\
&\leq
\sum_{k=1}^\infty r^k\left\|[A_{(\alpha)}-B_{(\alpha)}:\ |\alpha|=k]\right\|\\
&\leq \sum_{k=1}^\infty r^k \|B_{(0)}-A_{(0)}\| \\
&\leq \|B_{(0)}-A_{(0)}\|\frac{r}{1-r} \\
&\leq \frac{1}{2}\|B_{(0)}-A_{(0)}\|
\end{split}
\end{equation*}
for any $r\in\left[0,\frac{1}{3}\right]$.
Therefore, the second inequality of the theorem is proved.

Now, notice that if $A_{(0)}=B_{(0)}$, then adding $\epsilon I$,\ $\epsilon>0$, to  inequality \eqref{poz2} and applying the first part of the proof, we deduce that
$$
          \sum_{|\alpha|=k} (B_{(\alpha)}- A_{(\alpha)})^*(B_{(\alpha)}- 
      A_{(\alpha)}) 
         \leq  \epsilon^2 I  
$$     
for any $\epsilon>0$. This implies $A_{(\alpha)}=B_{(\alpha)}$, for any $\alpha\in \FF_n^+$. The proof is complete.
\end{proof}

We remark  that Theorem \ref{f<g} remains true if  the conditions  on the coefficients \linebreak $\cA:=\{A_{(\alpha)}\}_{\alpha\in \FF_n^+}$   and $\cB:=\{B_{(\alpha)}\}_{\alpha\in \FF_n^+}$  are  replaced by the following weaker conditions:
\begin{enumerate}
\item[(i)] $\limsup\limits_{k\to \infty}\left(\sum\limits_{|\alpha|=k} \|A_{(\alpha)}h\|^2\right)^{1/2k}\leq 1$ \ for any $h\in \cH$, and a similar inequality for $\cB$;
\item[(ii)] $\left< H_\cA(rS_1,\ldots, rS_n)p,p\right>\leq
\left< H_\cB(rS_1,\ldots, rS_n)p,p\right> $ for any $r\in [0,1)$ and  $p\in \cP\otimes \cH$.
\end{enumerate}

Now, we can prove the following result for  selfadjoint harmonic polynomials.
 
\begin{theorem}\label{harm-joint}
Let $m=2,3,\ldots,$ and  let   $\cA:=\{A_\alpha\}_{|\alpha|\leq m-1}$ and 
$\cB:=\{B_\alpha\}_{|\alpha|\leq m-1}$ be sequences of operators  such that $A_{(0)}=A_{(0)}^*$, $B_{(0)}=B_{(0)}^*$, and 
$$
H_\cA(S_1,\ldots, S_n)\leq H_\cB(S_1,\ldots, S_n).
$$
Then
\begin{equation}\label{Wi}
w\left( B_{(\alpha)}^*-A_{(\alpha)}^*:\ |\alpha|=k\right)\leq 
\|B_{(0)}-A_{(0)}\|\cos\frac {\pi} {\left[\frac{m-1} {k}\right]+2}
\end{equation}
 for  $1\leq k\leq m-1$,  where $[x]$ is the integer part of $x$, and
\begin{equation*}
\sum_{k=1}^{m-1}r^k w\left(A_{(\alpha)}^*:\ |\alpha|=k\right)\leq \frac{1}{2} \|B_{(0)}-A_{(0)}\|+
\sum_{k=1}^{m-1}r^k w\left(B_{(\alpha)}^*:\ |\alpha|=k\right)
\end{equation*}
for any 
$r\in [0,t_m]$, where $t_m$ is the solution of the equation
\eqref{eq}.
\end{theorem}
\begin{proof}
Since
\begin{equation*}  \sum_{k=1}^{m-1}\sum_{|\alpha|=k} r^{|\alpha|} S_\alpha^*\otimes 
      \left (B_{(\alpha)}^*- A_{(\alpha)}^*\right) 
       + I\otimes \left(B_{(0)}-A_{(0)}\right)+ 
       \sum_{k=1}^{m-1} \sum_{|\alpha|=k} r^{|\alpha|} S_\alpha\otimes 
      \left( B_{(\alpha)}- A_{(\alpha)}\right)\geq 0         \end{equation*} 
for $r\in [0,1)$,  inequality \eqref{joint-ine} implies
\eqref{Wi}.      
 Since the joint numerical radius is a norm and using the  inequality \eqref{Wi},  we have
   \begin{equation*}\begin{split}
\sum_{k=1}^{m-1}r^k w\left(A_{(\alpha)}^*:\ |\alpha|=k\right)&
-\sum_{k=1}^{m-1}r^k w\left(B_{(\alpha)}^*:\ |\alpha|=k\right)\\
&\leq \sum_{k=1}^{m-1}r^k w\left(B_{(\alpha)}^*-A_{(\alpha)}^*:\ |\alpha|=k\right)\\
&\leq 
 \|B_{(0)}-A_{(0)}\|\sum_{k=1}^{m-1} r^k\cos\frac{\pi} {\left[\frac{m-1}{k}\right]+2}\\
&\leq \|B_{(0)}-A_{(0)}\|\sum_{k=1}^{m-1} t_m^k\cos\frac{\pi} {\left[\frac{m-1}{k}\right]+2}\\
&\leq  \frac{1}{2} \|B_{(0)}-A_{(0)}\|
  \end{split}
\end{equation*}
for any $r\in [0, t_m]$.
This completes the proof.
\end{proof}

We remark that  one can obtain  versions of Theorem \ref{f<g}  and Theorem \ref{harm-joint} for the spatial tensor products $C_{red}^*(\FF_n)\otimes B(\cH)$ and $C^*(\FF_n)\otimes B(\cH)$. In this way we can  provide operator-valued generalizations of Corollary \ref{CUU} and Corollary \ref{CUU3}.

Given $m=2,3,\ldots, \infty$, and a sequence $\cA:=\{a_\alpha\}_{|\alpha|\leq m-1}$ of complex numbers with \linebreak
$a_0=\bar a_0$, define
$$H_\cA(S_1,\ldots, S_n):=
\sum_{1\leq |\alpha|\leq m-1} \bar a_\alpha S_\alpha^*+
a_0I+\sum_{1\leq |\alpha|\leq m-1} a_\alpha S_\alpha.
$$
When $m=\infty$, we assume that the series $\sum\limits_{k=1}^\infty\left\|\sum\limits_{|\alpha|=k}  r^{|\alpha|}a_\alpha S_\alpha\right\|$  is  convergent for any $r\in [0,1)$.
The following result is a scalar version of Theorem \ref{f<g} and Theorem \ref{harm-joint}.

\begin{corollary}
Let $m=2,3,\ldots, \infty$ and  let $\cA:=\{a_\alpha\}_{|\alpha|\leq m-1} $
and $\cB:=\{b_\alpha\}_{|\alpha|\leq m-1} $ be sequences 
of complex numbers such that $a_0=\bar a_0$, $b_0=\bar b_0$, and 
$$
H_\cA(S_1,\ldots, S_n)\leq 
H_\cB(S_1,\ldots, S_n) \quad \text{ if } \ m<\infty,
$$
and, if $ m=\infty$,
 
$$
H_\cA(rS_1,\ldots, rS_n)\leq 
H_\cB(rS_1,\ldots, rS_n) \quad \text{ for any } r\in [0,1). $$

Then
\begin{equation*}
\left(\sum_{|\alpha|=k}|b_\alpha-a_\alpha|^2\right)^{1/2}\leq 
(b_0-a_0)\cos\frac{\pi} {\left[\frac{m-1}{k}\right]+2}
\end{equation*}
for $1\leq k\leq m-1$, where $[x]$ is the integer part of $x$, and 
\begin{equation*}
\frac{a_0}{2}+\sum_{k=1}^{m-1}r^k\left(\sum_{|\alpha=k} |a_\alpha|^2\right)^{1/2}\leq 
\frac{b_0}{2}+\sum_{k=1}^{m-1}r^k \left(\sum_{|\alpha|=k}  |b_\alpha|^2\right)^{1/2}
\end{equation*}
for any $r\in[0,t_m]$, where $t_m$ is the solution of the equation \eqref{eq} if $m<\infty$ and $t_\infty:=\frac{1}{3}$.
\end{corollary}

Now, we can deduce Bohr inequalities for two real-valued functions $f,g\in L^1(\TT)$ satisfying the inequality $f\leq g$.

\begin{corollary} Let $f(e^{it})=\sum\limits_{k=-m+1}^{m-1} a_k e^{ikt}$ and $g(e^{it})=\sum\limits_{k=-m+1}^{m-1} a_k e^{ikt}$ be real-valued functions in $L^1(\TT)$, 
where  $m=2,3,\ldots, \infty$. 
If $f\leq g$, then
$$
 \frac{a_0} {2}+\sum_{k=1}^{m-1} r^k{|a_k|} \leq
\frac{b_0} {2}+\sum_{k=1}^{m-1} r^k{|b_k|}
$$
for any $r\in [0, t_m]$, where $t_m$ is the solution of the equation \eqref{eq} if $m<\infty$  and $t_\infty:=\frac{1}{3}$.
\end{corollary}

\begin{proof} Taking the harmonic extension of $f$ and $g$ in the unit disc, we obtain
$f(re^{it})\leq g(re^{it})$ for $0\leq r<1$ and 
$t\in [-\pi, \pi]$.
 Notice  that, for each $m=2,3,\ldots, \infty$,  $$H_{\cA}(rS):=\sum_{k=1}^{m-1} r^k \bar a_k {S^*}^k +a_0I+
\sum_{k=1}^{m-1} r^k  a_k S^k,\quad r\in [0,1),
$$
is in the $C^*$-algebra generated by the unilateral shift $S$, acting 
  on $H^2(\TT)$.
Similarly, we can define $H_{\cB}(rS)$, where $\cB:=\{b_k\}_{k=0}^{m-1}$.
Now, notice that, for any $h(e^{it})  \in H^2(\TT)$, we have 
\begin{equation*}
\begin{split}
\left< \left(H_{\cB}(rS)-H_{\cA}(rS) \right)h(e^{it}),h(e^{it})\right>_{H^2(\TT)}
&= 
\left< \left(g(re^{it})-f(re^{it})  \right)h(e^{it}),h(e^{it})\right>_{H^2(\TT)}\\
&=
\frac{1} {2\pi}\int_{-\pi}^\pi (g(re^{it})-f(re^{it})|h(e^{it})|^2dt\geq 0. 
\end{split}
\end{equation*}
Therefore, $H_{\cB}(rS)\geq H_{\cA}(rS)$ for  $0\leq r<1$.
Applying now Theorem \ref{f<g} (when $m=\infty$)  and Theorem \ref{harm-joint} (when $m=2,3,\ldots $), the result follows.
\end{proof}

\begin{corollary} Let $f:=\sum\limits_{k=-m+1}^{m-1} a_ke^{ikt}$ and $g:=\sum\limits_{k=-m+1}^{m-1} b_ke^{ikt}$ be  real-valued functions in $L^1(\TT)$, where  $m=2,3,\ldots, \infty$.
If there exists $r\in [0,t_m]$, where $t_m$ is the solution of the equation \eqref{eq} if $m<\infty$ and $t_\infty:=\frac{1}{3}$, such that 
\begin{equation*}
\frac{a_0}{2}+\sum_{k=1}^{m-1} r^{k} |a_k|>
\frac{b_0}{2}+\sum_{k=1}^{m-1} r^{k} |b_k|,
\end{equation*}
then $f>g$ on a set of positive measure.
\end{corollary}

We mention that,  when  $m=\infty$,  the last two corollaries    were obtained by Paulsen and Singh \cite{PS}.

\bigskip

\section{Bohr inequalities for analytic functions on the unit ball of $\CC^n$}

 A    multivariable commutative analogue 
         of the Hardy space $H^\infty(\DD)$  is the algebra \linebreak
         $W_n^\infty:=P_{F_s^2(H_n)} F_n^\infty|_{F_s^2(H_n)}$, 
        the compression  of $ F_n^\infty$ to the symmetric Fock space 
         ${F_s^2(H_n)} \subset F^2(H_n)$.       
 We proved in  \cite{APo}  that $W_n^\infty$ is the $w^*$-closed algebra generated by  the creation
      operators 
         $B_i:=P_{F_s^2(H_n)} S_i|_{F_s^2(H_n)}, \ i=1,\dots, n$, and the 
      identity. 
         Arveson \cite{Arv} showed    
         that the algebra $W_n^\infty$ can be seen as the   multiplier algebra 
      of the 
      reproducing kernel Hilbert space 
      with reproducing kernel $K_n: \BB_n\times \BB_n\to \CC$ defined by 
       $$ 
       K_n(z,w):= {\frac {1} 
      {1-\langle z, w\rangle_{\CC^n}}}, \qquad z,w\in \BB_n, 
      $$ 
      where $\BB_n$ is the open unit ball of $\CC^n$.

Let ${\bf p}:=(p_1,\ldots, p_n)$ be a multi-index in $\ZZ_+^n$. We denote 
$|{\bf p}|:=p_1+\cdots + p_n$ and ${\bf p} !:={ p}_1 !\cdots { p}_n !$.
If $\lambda:=(\lambda_1,\ldots,\lambda_n)$, then we set 
$\lambda^{\bf p}:=\lambda_1^{p_1}\cdots \lambda_n^{p_n}$ and 
$$
(\lambda^{\bf p})_{\text{\rm sym}} (S_1,\ldots, S_n):=\frac {{\bf p}!} {|{\bf p}|! }\sum_{\alpha\in \Lambda_{\bf p}} S_\alpha,
$$
where
$$
\Lambda_{\bf p}:=\{\alpha\in \FF_n^+: \lambda_\alpha= \lambda^{\bf p} \text{ for any } \lambda\in \BB_n\}
$$
and $S_1,\ldots, S_n$ are the left creation operators on the Fock space $F^2(H_n)$.
Notice that card\,$\Lambda_{\bf p}=\frac {|{\bf p}|!} 
{{\bf p}!}$.
If $p(\lambda):=\sum \lambda^{\bf p} A_{({\bf p})}$, 
$A_{({\bf p})}\in B(\cH)$,  is an operator-valued polynomial on $\BB_n$, the open unit ball
of $\CC^n$,  then the
symmetrized functional calculus 
$p_{\text{\rm sym}}(S_1,\ldots, S_n)\in F_n^\infty\bar{\otimes} B(\cH)$
is defined by
$$
p_{\text{\rm sym}}(S_1,\ldots, S_n):=\sum [(\lambda^{\bf p})_{\text{\rm sym}} (S_1,\ldots, S_n)]\otimes A_{({\bf p})}.
$$

  In this section we obtain  Wiener and Bohr   inequalities for the algebra $W_n^\infty$  and a class of operator-valued analytic functions on the open unit ball of $\CC^n$.  
 
\begin{theorem}\label{commutative}
Let $f(\lambda_1,\ldots,\lambda_n):=\sum\limits_{\bf p\in \ZZ_+^n}
\lambda^{\bf p} A_{(\bf p)}, \ A_{(\bf p)}\in B(\cH)$,
be  an operator-valued analytic  function on $\BB_n$ 
such that
\begin{equation}\label{sup}
\limsup_{k\to \infty}\left\|
\sum\limits_{{\bf p}\in \ZZ_+^n,|{\bf p}|=k}
 \frac 
{|{\bf p}|!}{{\bf p}!} 
A^*_{({\bf p})}A_{({\bf p})}\right\|^{1/2k}\leq 1,
\end{equation}
 $f(0)\geq 0$,  and
$$
  \text{\rm Re} 
f_{\text{\rm sym}}(rS_1,\ldots, rS_n)\leq I\quad \text{ for } \  0\leq r<1.
$$
Then the following statements hold.
\begin{enumerate}
\item[(i)]
$\left\| \sum\limits_{{\bf p}\in \ZZ_+^n, |{\bf p}|=k} \frac {{\bf p}!}{|{\bf p}|!} 
A^*_{({\bf p})}A_{({\bf p})}\right\|^{1/2}
\leq 2\|I-A_{(0)}\|$ \ for $k=1,2,\ldots$.
\item[(ii)]
$\sum\limits_{k=0}^\infty r^k 
\left\| \sum\limits_{{\bf p}\in \ZZ_+^n,|{\bf p}|=k} \frac 
{{\bf p}!}{|{\bf p}|!} 
A^*_{({\bf p})}A_{({\bf p})}\right\|^{1/2}
\leq \|A_{(0)}\|+\|I-A_{(0)}\|$ \ 
for $0\leq r\leq \frac {1} {3}$.
\item[(iii)]
If $\{Y_{({\bf p})}\}_{{\bf p}\in \ZZ_+^n}$ is a sequence of operators in $B(\cK)$ such that
$$
\|Y_{(0)}\|\leq 1 \ \text{ and } \ 
\sum\limits_{k=0}^\infty   
\left\| \sum\limits_{{\bf p}\in \ZZ_+^n,|{\bf p}|=k} \frac {|{\bf p}|!} 
{{\bf p}!}Y^*_{(\bf p)}Y_{(\bf p)}\right\|^{1/2}
\leq \frac {1} {2},
$$
then $$\left\|\sum\limits_{k=0}^\infty   
 \left(\sum\limits_{{\bf p}\in \ZZ_+^n,|{\bf p}|=k} B^{\bf p}\otimes 
A_{({\bf p})}\otimes Y_{({\bf p})}\right)\right\|\leq 1,
$$
where $B_1,\ldots, B_n$ are the creation operators on the symmetric Fock space.
In particular, 
$$\left\|\sum\limits_{k=0}^\infty   
\left( \sum\limits_{{\bf p}\in \ZZ_+^n, |{\bf p}|=k} \lambda^{{\bf p}} 
A_{({\bf p})}\otimes Y_{(\bf p)}\right)\right\|\leq 1
 \ \text{ for any }\   \lambda\in \BB_n.$$
\item[(iv)] If $Y_1,\ldots, Y_n\in B(\cK)$ are commuting and 
$\|[Y_1^*,\ldots, Y_n^*]\|\leq \frac {1} {3}$, then
$$\left\|\sum\limits_{k=0}^\infty   
 \left(\sum\limits_{{\bf p}\in \ZZ_+^n,|{\bf p}|=k} B^{{\bf p}}\otimes Y^{{\bf p}}\otimes 
A_{({\bf p})}\right)\right\|\leq 1.
$$
\item[(v)] If $(T_1,\ldots, T_n)\in [B(\cK)^n]_{1/3}$ and $T_1,\ldots, T_n$ are commuting operators, then
$$\sum\limits_{k=0}^\infty \left\|  
 \sum\limits_{{\bf p}\in \ZZ_+^n, |{\bf p}|=k} T^{\bf p}\otimes 
A_{({\bf p})} \right\|\leq \|A_{(0)}\|+\|I-A_{(0)}\|.
$$
\end{enumerate}
\end{theorem}

\begin{proof} Notice that 
$$f_{\text{\rm sym}}(rS_1,\ldots, rS_n)=\sum_{k=0}^\infty \sum_{|\alpha|=k} r^{|\alpha|} S_\alpha \otimes  C_{(\alpha)},
$$
where $C_{(0)}:=A_{(0)}$   and $C_{(\alpha)}:= 
 \frac {{\bf p}!}{|{\bf p}|!}A_{({\bf p})}$ for ${\bf p}\in \ZZ_+^n$, ${\bf p}\neq (0,\ldots, 0)$, and $\alpha\in \Lambda_{\bf p}$.
We remark that condition \eqref{sup} implies that $f_{\text{\rm sym}}(rS_1,\ldots, rS_n)$ is norm convergent for each $r\in [0,1)$.
It is clear that, for each $k=1,2,\ldots, $ we have
\begin{equation*}
\begin{split}
\sum_{|\alpha|=k}C_{(\alpha)}^*C_{(\alpha)}&=
\sum\limits_{{\bf p}\in \ZZ_+^n,|{\bf p}|=k}\left(\sum_{\alpha\in \Lambda_{\bf p}} C_{(\alpha)}^*C_{(\alpha)}\right)\\
&=\sum\limits_{{\bf p}\in \ZZ_+^n,|{\bf p}|=k}
 \frac 
{{\bf p}!}{|{\bf p}|!} 
A^*_{({\bf p})}A_{({\bf p})}.
 \end{split}
\end{equation*}
Apply now Theorem \ref{bohr-gen} to the  holomorphic function 
$F(X_1,\ldots, X_n):=f_{\text{\rm sym}}(X_1,\ldots, X_n)$,\ $(X_1,\ldots, X_n)\in [B(\cX)^n]_{<1}$, and notice that part (i) and (ii) follow from part (ii) and part (vi) of 
Theorem \ref{bohr-gen}.
To prove part (iii), define 
the sequence  $\{Z_{(\alpha)}\}_{\alpha\in \FF_n^+}\subset B(\cK)$ by $Z_{(0)}:=Y_{(0)}$ and $Z_{(\alpha)}:=Y_{({\bf p})}$, where 
${\bf p}\in \ZZ_+^n$, ${\bf p}\neq (0,\ldots, 0)$, and $\alpha\in \Lambda_{\bf p}$.
Notice also that 
$$
\sum_{k=0}^\infty \left\| \sum_{|\alpha|=k}Z_{(\alpha)}^*Z_{(\alpha)} \right\|^{1/2}=
\sum\limits_{k=0}^\infty   
\left\| \sum\limits_{{\bf p}\in \ZZ_+^n,|{\bf p}|=k} \frac {|{\bf p}|!} 
{{\bf p}!}Y^*_{(\bf p)}Y_{(\bf p)}\right\|^{1/2}
\leq \frac {1} {2},
$$
Using  Theorem \ref{noncom} part (ii) when $F(X_1,\ldots, X_n)=f_{\text{\rm sym}}(X_1,\ldots, X_n)$, we deduce that 
$$
\left\| \sum_{k=0}^\infty \sum_{|\alpha|=k} S_\alpha\otimes C_{(\alpha)}\otimes Z_{(\alpha)}\right\|\leq 1,
$$
where the series converges in the norm topology.
Taking the compression to the symmetric Fock space  and noticing that 
\begin{equation*}
\begin{split}
\sum_{k=0}^\infty \sum_{|\alpha|=k} B_\alpha\otimes C_{(\alpha)}\otimes Z_{(\alpha)}&=
\sum\limits_{{\bf p}\in \ZZ_+^n,|{\bf p}|=k}\left(\sum_{\alpha\in \Lambda_{\bf p}}\frac {{\bf p}!}{|{\bf p}|!}B^{\bf p} \otimes A_{\bf p} \otimes Y_{\bf p}\right)\\
&=
\sum\limits_{{\bf p}\in \ZZ_+^n,|{\bf p}|=k} B^{\bf p} \otimes A_{\bf p} \otimes Y_{\bf p},
\end{split}
\end{equation*}
we conclude the proof of part (iii) of the theorem.

Now, applying Theorem \ref{bohr-gen} part (i) to 
$f_{\text{\rm sym}}(X_1,\ldots, X_n)$, we obtain
$$  
\left\|\sum_{k=0}^\infty \sum_{|\alpha|=k} S_\alpha \otimes Y_\alpha\otimes C_{(\alpha)}\right\|\leq 1.
$$
Taking again the compression to the symmetric Fock space, we deduce (iv).
To prove part (v), we apply Theorem \ref{Bo-ge}
to $f_{\text{\rm sym}}(X_1,\ldots, X_n)$ and deduce  the inequality
$$
\sum_{k=0}^\infty\left\| \sum_{|\alpha|=k} T_\alpha\otimes C_{(\alpha)}\right\|\leq \|C_{(0)}\|+\|I-C_{(0)}\|.
$$
On the other hand, since $T_1,\ldots, T_n$ are commuting,  we deduce that
$$
\sum_{k=0}^\infty\left\| \sum_{|\alpha|=k} T_\alpha\otimes C_{(\alpha)}\right\|=
\sum\limits_{k=0}^\infty \left\|  
 \sum\limits_{{\bf p}\in \ZZ_+^n, |{\bf p}|=k} T^{\bf p}\otimes 
A_{({\bf p})} \right\|.
$$
This completes the proof.
\end{proof}

\begin{corollary}\label{commu}
Let $f(\lambda_1,\ldots,\lambda_n):=\sum\limits_{{\bf p}\in \ZZ_+^n}
\lambda^{{\bf p}} a_{{\bf p}}$, \ $ a_{{\bf p}}\in \CC$,
be  an analytic function on $\BB_n$  such that
 condition \eqref{sup} holds,\ 
 $f(0)\geq 0$, and  
    $$
  \text{\rm Re} 
f_{\text{\rm sym}}(rS_1,\ldots, rS_n)\leq I\ \  \text{ for } 0\leq r<1.
$$
 Then
$$\sum\limits_{k=0}^\infty \left|  
 \sum\limits_{{\bf p}\in \ZZ_+^n, |{\bf p}|=k} \lambda^{\bf p} |a_{{\bf p}}| \right|\leq  1
$$
for any $\lambda:=(\lambda_1,\ldots, \lambda_n)\in \BB_n$ with
$\|\lambda\|_2\leq \frac {1} {3}$.
\end{corollary}

Theorem \ref{pol} can be used to obtain a version of Theorem \ref{commutative} for operator-valued polynomias in $\BB_n$.
Here is    a scalar version.

\begin{corollary}
Let $m=2,3,\ldots$, and let
$$p(\lambda_1,\ldots,\lambda_n):=\sum\limits_{{\bf p}\in \ZZ_+^n, |{\bf p}|\leq m-1}
\lambda^{{\bf p}} a_{{\bf p}}, \quad a_{{\bf p}}\in \CC,
$$ 
 be  an  analytic polynomial on $\BB_n$,  such that
 $p(0)\geq 0$ and  
    $
 \  \text{\rm Re}\, 
p_{\text{\rm sym}}(S_1,\ldots, S_n)\leq I $. 
 Then 
$$\sum\limits_{k=0}^{m-1} \left|  
 \sum\limits_{{\bf p}\in \ZZ_+^n, |{\bf p}|=k} \lambda^{\bf p} |a_{{\bf p}}| \right|\leq  1
$$
for any $\lambda:=(\lambda_1,\ldots, \lambda_n)\in \BB_n$ with
$\|\lambda\|_2\leq t_m$, where $t_m$ is the solution of the equation \eqref{eq}.
\end{corollary}

Denote by $H_{\text{\rm sym}}^\infty(\BB_n)$ the set of all 
analytic functions on $\BB_n$ with scalar coefficients
 $$
f(\lambda_1,\ldots,\lambda_n):=\sum\limits_{\bf p\in \ZZ_+^n}
\lambda^{\bf p} a_{\bf p}, \quad  a_{\bf p}\in \CC,
$$
such that
\begin{equation*}
\limsup_{k\to \infty}\left(
\sum\limits_{{\bf p}\in \ZZ_+^n,|{\bf p}|=k}
 \frac 
{|{\bf p}|!}{{\bf p}!} 
|a_{\bf p}|^2\right)^{1/2k}\leq 1
\end{equation*}
and
$$\|f\|_{\text{\rm sym}}:=\sup_{0\leq r<1}\left\| f_{\text{\rm sym}}(rS_1,\ldots, rS_n)\right\|<\infty.
$$

According to \cite{Po-holomorphic},   $H_{\text{\rm sym}}^\infty(\BB_n)$ is a Banach space with respect to the the norm
$\|\cdot\|_{\text{\rm sym}}$, which contains all analytic polynomials on $\BB_n$.

\begin{corollary}
Let $m=2,3,\ldots, \infty$, and let
$$f(\lambda_1,\ldots,\lambda_n):=\sum\limits_{{\bf p}\in \ZZ_+^n, |{\bf p}|\leq m-1}
\lambda^{{\bf p}} a_{{\bf p}}, \quad a_{{\bf p}}\in \CC,
$$ 
 be  an analytic function in $H_{\text{\rm sym}}^\infty(\BB_n)$.    Then 
$$\sum\limits_{k=0}^{m-1} \left|  
 \sum\limits_{{\bf p}\in \ZZ_+^n, |{\bf p}|=k} \lambda^{\bf p} |a_{{\bf p}}| \right|\leq  \|f\|_{\text{\rm sym}}
$$
for any $\lambda:=(\lambda_1,\ldots, \lambda_n)\in \BB_n$ with
$\|\lambda\|_2\leq t_m$, where $t_m$ is the solution of the equation \eqref{eq} if $m<\infty$ and $t_\infty=\frac{1}{3}$.
\end{corollary}

       %

      \end{document}